\documentclass[11pt]{article}
\usepackage{latexsym,amssymb}
\usepackage{psfrag}
\usepackage{amsmath}
\usepackage{color}
\newtheorem{Theorem}{\bf Theorem}[section]
\newtheorem{Lemma}{\bf Lemma}[section]
\newtheorem{Proposition}{\bf Proposition}[section]
\newtheorem{Corollary}{\bf Corollary}[section]
\newtheorem{Remark}{\bf Remark}[section]
\newtheorem{Example}{\bf Example}[section]
\newtheorem{Definition}{\bf Definition}[section]

\newenvironment{theorem}{\begin{Theorem}$\!\!\!$}{\end{Theorem}}
\newenvironment{lemma}{\begin{Lemma}$\!\!\!$}{\end{Lemma}}
\newenvironment{proposition}{\begin{Proposition}$\!\!\!$}{\end{Proposition}}

\numberwithin{equation}{section}
\begin{document}

\title{The heat kernel of a Schr\"odinger  operator\\ 
with  inverse square  potential}
\author{Kazuhiro Ishige, Yoshitsugu Kabeya and El Maati Ouhabaz}
\date{}
\maketitle
\begin{abstract}{We consider the Schr\"odinger operator $H = -\Delta + V(|x|)$ with radial potential $V$ which may have singularity at $0$ and a quadratic decay at infinity. First, we study the structure of positive harmonic functions of $H$ and give their precise behavior. Second, under quite general conditions we prove an upper  bound for the correspond heat kernel $p(x,y,t)$ of the type
\begin{equation*}
0 < p(x,y,t)
\le C\,t^{-\frac{N}{2}}
\frac{U(\min\{|x|,\sqrt{t}\})U(\min\{|y|,\sqrt{t}\})}{U(\sqrt{t})^2}
\exp\left(-\frac{|x-y|^2}{Ct}\right)
\end{equation*}
for all $x$, $y\in{\bf R}^N$ and $t>0$, where $U$ is a positive harmonic function of $H$. Third, if $U^2$ is an $A_2$ weight on ${\bf R}^N$, then we  prove a lower bound of a similar type.}
\end{abstract}
\section{Introduction}\label{section:1}

Heat kernel bounds of differential  operators on domains of ${\bf R}^N$ or Riemannian  manifolds have attracted attention in recent years. We refer the reader for an account on this to the monographs of Davies \cite{Dav}, Grigor'yan \cite{Gri} and Ouhabaz \cite{Ouh}.
Typically, for a second order  differential  elliptic operator $H$, the associated heat kernel $p(x,y,t)$ (i.e. the integral kernel of the semigroup $e^{-tH}$ generated by $-H$, or the fundamental solution to the heat equation associated with $H$)  satisfies  in many cases the following upper bound
\begin{equation*}
| p(x,y,t) | \le \frac{C}{\sqrt{|B(x,\sqrt{t})|} \sqrt{|B(y, \sqrt{t})|}} e^{- \frac{| x- y|^2}{ct}}, 
\end{equation*}
where $|B(x,r)|$ denotes the volume of the open ball of the manifold with center $x$ and radius $r$ and $|x - y|$ denotes the Riemannian distance between the two points $x$ and $y$. In the Euclidean setting (i.e.  ${\bf R}^N$) the above estimate reduces to
\begin{equation*}
| p(x,y,t) | \le C t^{-N/2} e^{- \frac{| x- y|^2}{ct}}.
\end{equation*}
These bounds are referred to as Gaussian upper bounds for $p(x,y,t)$. Such bounds have been studied in many situations. They play an important role in several  problems. For example, they are used in harmonic analysis in order to prove boundedness of some singular integral operators such as Riesz transforms or spectral multipliers, in spectral theory in order to prove $p$-independence of the spectrum, to prove maximal regularity for the evolution equation, and so on. 
For all this we refer to Chapter 7 in \cite{Ouh} and references there.

There are however many cases where such upper bound cannot hold. 
A typical and important example is the Schr\"odinger operator with inverse square potential, i.e., 
$$H = -\Delta+\frac{\lambda}{|x|^2}, $$
where $-(N-2)^2/4\le\lambda<0$. 
It is well known that the semigroup $e^{-tH}$ does not act on $L^p({\bf R}^N)$ for 
$p$ outside a certain symmetric interval around $2$ whose length depends on the constant $\lambda$. See Liskevich, Sobol and Vogt \cite{LSV}. 
Therefore, the corresponding heat kernel $p(x,y,t)$ does not satisfy the above classical Gaussian bound.  It was proved by Milman and Semenov \cite{MS2},  and later by Liskevich and Sobol \cite{VS} that the heat kernel satisfies
$$
0 < p(x,y,t)
\le C\,t^{-\frac{N}{2} + \sigma}
(\min\{|x|,\sqrt{t}\})^{-\sigma}(\min\{|y|,\sqrt{t}\})^{-\sigma}
\exp\left(-\frac{|x-y|^2}{Ct}\right)
$$
for all $x$, $y\in{\bf R}^N$ and $t>0$, where 
$$
\sigma = \frac{N-2}{2} - \frac{1}{2} \sqrt{(N-2)^2 +4\lambda}.
$$  
See also Barbatis, Filippas and Tertikas  \cite{Bar}. The result in \cite{VS} deals with a more general class of operators in the sense that $\Delta$ is replaced by a divergence form operator with appropriate behavior of the coefficients. A lower bound of the same type was also proved in \cite{MS1} and \cite{MS2}. 
We observe that this upper bound can be rephrased as
\begin{equation}\label{eq:1.1}
0 < p(x,y,t)
\le C\,t^{-\frac{N}{2}}
\frac{U(\min\{|x|,\sqrt{t}\})U(\min\{|y|,\sqrt{t}\})}{U(\sqrt{t})^2}
\exp\left(-\frac{|x-y|^2}{Ct}\right), 
\end{equation}
where $U(x) = |x|^{-\sigma} $ and it turns out that $U$ is a positive harmonic function of $H$.  

Our aim in this paper is to prove the  bounds as in \eqref{eq:1.1} for a wide class of potentials. 
Thus we are led to consider first existence and behavior of positive harmonic functions. 

The behavior of positive harmonic functions for Schr\"odinger operators have been studied by Murata \cite{M01}. 
 He studied the structure of all positive harmonic functions for the elliptic operator $-\Delta+V(x)$
in the case where $V\in L^p_{{\rm loc}}({\bf R}^N)$ with some $p>N/2$ if $N\ge 2$ and $p>1$ if $N=1$. 
Furthermore, he classified the behavior of positive harmonic functions, 
in particular, in the case where $V$ is a radially symmetric function on ${\bf R}^N$ 
(see \cite[Section~3]{M01}). See also Remark~\ref{Remark:1.1}.

In the present paper we consider a more general class of possibly negative potentials. We assume that $N \ge 2$ and the radial potential $V$ is continuous on $(0, \infty)$ and satisfies
\begin{equation}
\label{eq:1.2}
\left\{
 \begin{array}{l}
 \displaystyle{\lim_{r\to 0}r^{-\theta}\left|r^2V(r)-\lambda_1\right|=0},
 \quad
 \displaystyle{\lim_{r\to\infty}r^\theta\left|r^2V(r)-\lambda_2\right|=0},\vspace{3pt}\\
 \mbox{where $\lambda_1$, $\lambda_2\in[\lambda_*,\infty)$ with $\lambda_*:=-(N-2)^2/4$,}
\end{array}
\right.
\end{equation} 
for some $\theta>0$. 
We also assume that 
the Schr\"odinger operator $H:=-\Delta+V$ is nonnegative, that is
$$
\int_{{\bf R}^N}\left[|\nabla\phi|^2+V\phi^2\right]\,dx\ge 0
\quad
\mbox{for all $\phi\in C^\infty_0({\bf R}^N\setminus\{0\})$}. 
$$
We first study the behavior of positive harmonic functions in the light of Murata's paper \cite{M01}.  The result will be then used to prove upper and lower estimate for the heat kernel $p(x,y,t)$. In order to state our  results we introduce some definitions and notation. 

We say that $H$ is subcritical 
if, for any $W\in C_0^\infty({\bf R}^N)$, 
$H-\epsilon W$ is nonnegative for  small enough $\epsilon>0$; 
$H$ is critical if $H$ is not subcritical. 
On the other hand, if $H$ is not nonnegative, then $H$ is said to be supercritical. 

For any $\lambda\in[\lambda_*,\infty)$, 
let $A^\pm(\lambda)$ be roots of the algebraic equation 
$\alpha^2+(N-2)\alpha-\lambda=0$
such that $A^-(\lambda)\le A^+(\lambda)$, that is 
\begin{equation}
\label{eq:1.3}
A^\pm(\lambda):=\frac{-(N-2)\pm\sqrt{D_\lambda}}{2},
\quad\mbox{where}\quad
D_\lambda:=(N-2)^2+4\lambda\ge 0. 
\end{equation}
Then $v(r) :=r^{A^\pm(\lambda)}$ satisfies 
$$
v''+\frac{N-1}{r}v'-\frac{\lambda}{r^2}v=0
\quad\mbox{in}\quad(0,\infty). 
$$
Furthermore, it follows that 
\begin{equation}
\label{eq:1.4}
A^-(\lambda)<-\frac{N-2}{2}<A^+(\lambda)\quad\mbox{if}\quad \lambda>\lambda_*,
\quad
A^\pm(\lambda)=-\frac{N-2}{2}\quad\mbox{if}\quad\lambda=\lambda_*. 
\end{equation}
%

For positive functions $f$ and $g$ defined on  $(0,R)$ for some $R>0$, 
we write 
$$
f(r)\thicksim g(r)\quad\mbox{as}\quad r\to 0
\quad\mbox{if}\quad
\lim_{r\to 0}\frac{f(r)}{g(r)}=1. 
$$ 
Similarly, for  positive functions $f$ and $g$ defined on  $(R,\infty)$ for some $R>0$, 
we write  
$$
f(r)\thicksim g(r)\quad\mbox{as}\quad r\to\infty
\quad\mbox{if}\quad
\lim_{r\to\infty}\frac{f(r)}{g(r)}=1. 
$$
Furthermore, for any two nonnegative functions 
$f_1$ and $f_2$ defined on a set $D$, 
we write
$$
f_1(r)\asymp f_2(r)\quad\mbox{for}\quad r\in D
$$ 
if there exists a positive constant $C$ such that 
$C^{-1}f_2(r)\le f_1(r)\le Cf_2(r)$ for all $r\in D$. 
\vspace{3pt}

Now we are ready to state the main results of this paper. 
The first theorem ensures the existence of positive harmonic functions 
for the operator $H=-\Delta+V$ and classifies the behavior of positive harmonic functions. 

\begin{theorem}
\label{Theorem:1.1} 
Let $N\ge 2$. 
Let $V$ be a continuous function on $(0,\infty)$ 
satisfying \eqref{eq:1.2}. 
\begin{itemize}
  \item[{\rm (1)}] 
  There exists a unique solution $U$ of 
  \begin{equation*}
  \tag{O}
  U''+\frac{N-1}{r}U'-V(r)\,U=0\quad\mbox{in}\quad(0,\infty),
  \end{equation*}
 with the property  $U(r)\thicksim r^{A^+(\lambda_1)}$ as $r\to 0$. 
  \item[{\rm (2)}] 
  For any solution $w$ of {\rm (O)} satisfying 
  $$
  w(r)=o(r^{A^-(\lambda_1)})\quad\mbox{if}\quad\lambda_1>\lambda_*,
  \qquad
  w(r)=o\left(r^{-\frac{N-2}{2}}|\log r|\right)\quad\mbox{if}\quad\lambda_1=\lambda_*,
  $$
  as $r\to 0$, 
  there exists a constant $c$ such that $w(r)=c\,U(r)$ on $(0,\infty)$, 
  where $U$ is as in ${\rm (1)}$. 
  \item[{\rm (3)}] 
  Assume that $H:=-\Delta+V$ is nonnegative. Then $U(r)>0$ on $(0,\infty)$ 
  and 
  \begin{itemize}
  \item[{\rm (a)}]
  $U(r)\thicksim c_*r^{A^+(\lambda_2)}$ as $r\to\infty$ if $H$ is subcritical and $\lambda_2>\lambda_*$,
  \item[{\rm (b)}]
  $U(r)\thicksim c_*r^{A^-(\lambda_2)}$ as $r\to\infty$ if $H$ is critical and $\lambda_2>\lambda_*$,
  \item[{\rm (c)}]
  $U(r)\thicksim c_*r^{-\frac{N-2}{2}}\log r$ as $r\to\infty$ if $H$ is subcritical and $\lambda_2=\lambda_*$,
  \item[{\rm (d)}]
  $U(r)\thicksim c_*r^{-\frac{N-2}{2}}$ as $r\to\infty$ if $H$ is critical and $\lambda_2=\lambda_*$,
  \end{itemize}
  for some $c_*>0$. 
  \item[{\rm (4)}] 
  Assume that $H$ is subcritical. Let $W\in C_0([0,\infty))$ be such that $W\ge 0$ and $W\not\equiv 0$ on $[0,\infty)$. 
  Set $H_\mu:=-\Delta+V-\mu W$ for $\mu\in{\bf R}$. 
  Then there exists $\mu_*>0$ such that 
  \begin{itemize}
  \item[{\rm (a)}] 
  $H_\mu$ is subcritical if $\mu<\mu_*$;
  \item[{\rm (b)}] 
  $H_\mu$ is critical if $\mu=\mu_*$;
  \item[{\rm (c)}] 
  $H_\mu$ is supercritical if $\mu>\mu_*$.
  \end{itemize}
\end{itemize}
\end{theorem}

\begin{Remark}
\label{Remark:1.1}
Let $N\ge 2$.
Let $V$ be a continuous function on $(0,\infty)$ 
satisfying \eqref{eq:1.2}. 
\begin{itemize}
\item[{\rm (i)}] 
In the case of $\lambda_1=0$,  
we see that $V\in L^{p/2}_{{\rm loc}}({\bf R}^N)$ for some $p>N/2$.  
Then Theorem~{\rm\ref{Theorem:1.1}} follows from Theorem~{\rm 3.1} in {\rm\cite{M01}}. 
\item[{\rm (ii)}] 
If $V(r)\ge\lambda r^{-2}$ on $[0,\infty)$ for some $\lambda>\lambda_*$, 
then $H$ is subcritical. This immediately follows from the Hardy inequality. 
\end{itemize}
\end{Remark}  

The next results are concerned with upper and lower bounds for the heat kernel $p(x,y,t)$  of $H = -\Delta + V$.  Recall that the heat kernel is the fundamental solution of
\begin{equation}
\label{eq:1.5}
\partial_t u=\Delta u-V(|x|)u
\quad\mbox{in}\quad{\bf R}^N\times(0,\infty).
\end{equation}
We prove  the following results.  
\begin{theorem}
\label{Theorem:1.2}
Let $N\ge 2$. 
Let $V$ be a continuous function on $(0,\infty)$ satisfying \eqref{eq:1.2} 
and $p=p(x,y,t)$ the fundamental solution of \eqref{eq:1.5}. 
Assume that $H:=-\Delta+V(|x|)$ is nonnegative 
and let $U$ be as in Theorem~{\rm\ref{Theorem:1.1}}.
If $\omega(x):= U(|x|)^2$ is an $A_2$~weight on ${\bf R}^N$, 
then there exist  positive constants $C_1$ and $C_2$ such that 
\begin{equation*}
\begin{split}
\frac{C_1 U(x) U(y)}{\sqrt{\omega(B(x,\sqrt{t}))}\sqrt{\omega(B(y,\sqrt{t}))}} & \exp\left(-\frac{|x-y|^2}{C_1t}\right) 
\le p(x,y,t)\le\\
 & \frac{C_2 U(x) U(y)}{\sqrt{\omega(B(x,\sqrt{t}))}\sqrt{\omega(B(y,\sqrt{t}))}}\exp\left(-\frac{|x-y|^2}{C_2t}\right)
\end{split}
\end{equation*}
for all $x$, $y\in{\bf R}^N$ and $t>0$. Here  
$\omega(B(x,\sqrt{t})):=\int_{B(x,\sqrt{t})}\omega(z)\,dz$. 
\end{theorem}

We shall see  in the proof that the upper bound can be made slightly more precise in the sense that the constant 
$C_2$ could chosen to be arbitrary close to $4$. Indeed we prove that 
\begin{equation*}
p(x,y,t)\le
\frac{C_\epsilon U(x)U(y)}{\sqrt{\omega(B(x,\sqrt{t}))}\sqrt{\omega(B(y,\sqrt{t}))}}\exp\left(-\frac{|x-y|^2}{(4+\epsilon)t}\right).
\end{equation*}
for every $\epsilon > 0$. The constant $C_\epsilon$ is independent of $x$, $y$ and $t$. 

Note that 
in Theorem~\ref{Theorem:1.2}, 
$\omega\equiv U^2$ is an $A_2$~weight on ${\bf R}^N$ if $A^+(\lambda_1)< N/2$ and 
$$
C^{-1}r^{-A_1}\le U(r)\le Cr^{-A_2},\qquad r\ge 1,
$$
for some $A_1$ and $A_2$ such that $-N/2<A_2\le A_1<N/2$. 
Next we weaken the $A_2$-assumption on $\omega$ and obtain an upper Gaussian estimate
for  $p=p(x,y,t)$. 
\begin{theorem}\label{Theorem:1.3} 
Let $N\ge 2$.
Let $V$ be a continuous function on $(0,\infty)$ satisfying \eqref{eq:1.2}. 
Assume that $H:=-\Delta+V(|x|)$ is nonnegative 
and let $U$ be as in Theorem~{\rm\ref{Theorem:1.1}}.
Furthermore, if $H$ is critical, then we assume that 
$$
A^-(\lambda_2)>-\frac{N}{2}. 
$$
Then there exists a positive constant $C$ such that  
\begin{equation}
\label{eq:1.6}
0<p(x,y,t)
\le C\,t^{-\frac{N}{2}}
\frac{U(\min\{|x|,\sqrt{t}\})U(\min\{|y|,\sqrt{t}\})}{U(\sqrt{t})^2}
\exp\left(-\frac{|x-y|^2}{Ct}\right)
\end{equation}
for all $x$, $y\in{\bf R}^N$ and $t>0$. 
\end{theorem}
For the proof of Theorem~\ref{Theorem:1.3}, 
we apply a refinement of the technique developed in \cite{IIY01, IIY02, IK01,IK04} 
and construct supersolutions of \eqref{eq:1.5}. 
Furthermore, we combine the comparison principle with the standard arguments as given, 
for example, \cite[Section~6]{S}, 
and prove Theorem~\ref{Theorem:1.3}. 

The final result is an observation that for a non-necessarily radial positive potential $V$, if one knows that there exists 
a harmonic function $U$ which behaves as a polynomial on the whole ${\bf R}^N$, 
then  the Gaussian upper bound holds. 
More precisely,
\begin{Proposition}
\label{Proposition:1.1}
Suppose that $V \ge 0$ and that $H$ has a harmonic function $U$ 
satisfying 
$$
C_0\, | x |^\alpha \le U(x) \le C_0'\, |x|^\alpha,
\qquad x\in{\bf R}^N,
$$
for some $\alpha \ge 0$ and $C_0$, $C_0' > 0$.
Then
\begin {equation*}
p(x,y,t) \le C\,t^{-\frac{N}{2}}
\frac{U(\min\{|x|,\sqrt{t}\})U(\min\{|y|,\sqrt{t}\})}{U(\sqrt{t})^2}
\exp\left(-\frac{|x-y|^2}{Ct}\right)
\end{equation*}
for all $x$, $y\in{\bf R}^N$ and $t>0$. 
\end{Proposition}

The proof of the latter result uses the standard Caffarelli-Kohn-Nirenberg inequalities. The idea is classical and we work on the weighted 
space $L^2({\bf R}^N, |x|^\alpha dx)$. Then the Sobolev inequality on this weighted space (which is the Caffarelli-Kohn-Nirenberg inequality) 
allows us to obtain an appropriate $L^2({\bf R}^N, |x|^\alpha dx)-L^\infty({\bf R}^N, |x|^\alpha dx)$ decay of the semigroup. 
The standard perturbation method allows then to convert this decay into a Gaussian bound. 
This reasoning has already appeared in \cite{Bar} in the context of the Schr\"odinger operator $-\Delta + \frac{\lambda}{|x|^2}$. 

Note that the above results extend the results  from the papers \cite{Bar}, \cite{VS}, \cite{MS1} and \cite{MS2}  
mentioned above which deal with the case where $V = \frac{\lambda}{|x|^2}$.  
\section{Preliminaries}\label{section:2}
In this section we  recall some properties 
for parabolic equations with $A_2$ weight. 
Throughout this section and in the rest of the paper, we denote by $C$  generic positive constants 
which may have different values even within the same line.

Let $\omega$ be a nonnegative measurable function on  a domain $\Omega\subset{\bf R}^N$. 
Suppose that  $\omega$ is  an $A_2$~weight on $\Omega$, that is 
$\omega$, $\omega^{-1}\in L^1_{{\rm loc}}(\Omega)$ and 
$$
[\omega](\Omega) := \sup\left\{\int_E \omega\,dz\,\int_E \omega^{-1}\,dz\biggr/\left(\int_E\,dz\right)^2\,:\,
\mbox{$E$ is a ball in $\Omega$}\right\}<\infty. 
$$
Then $\omega(z)\, dz$ is a measure on $\Omega$ with the doubling property, that is 
\begin{equation}
\label{eq:2.1}
\omega(B(x,2r))\le C\omega(B(x,r))
\end{equation}
holds for all $x\in{\bf R}^N$ and $r>0$, where 
$\omega(B(x,r)):=\int_{B(x,r)}\omega(z)\,dz$. 
For further details on $A_2$~weights, see e.g., \cite{Stein}.

We denote by $L^p(\Omega,\omega\,dx)$ $(1\le p<\infty)$ 
the usual Lebesgue  spaces  with norm 
$$
\|f\|_{p,\omega;\,\Omega}:=
\left(\int_\Omega |f(z)|^p\omega(z)\,dz\right)^{\frac{1}{p}}.
$$
By  $H^1(\Omega,\omega\,dx)$ we denote the Sobolev space defined as  the completion of $C^\infty(\overline{\Omega})$ 
with respect to  the norm  
$$
\left(\int_\Omega(|f(z)|^2+|\nabla f(z)|^2)\omega(z)\,dz\right)^{\frac{1}{2}}. 
$$

Consider the degenerate parabolic equation 
\begin{equation}
\label{eq:2.2}
\partial_t v=\frac{1}{\omega}\mbox{div}\,(\omega\nabla v)+cv\quad\mbox{in}\quad\Omega\times I,
\end{equation}
where $I$ is an open interval of ${\bf R}$ and $c\in L^\infty(I:L^\infty(\Omega))$. 
We say that a measurable function $v$ on $\Omega\times I$ is a solution of \eqref{eq:2.2} if 
$$
v\in L^\infty(I:L^2(\Omega),w\,dz))\,\cap\,L^2(I:H^1(\Omega,w\,dz)) 
$$
and $v$ satisfies 
$$
\int_I\int_\Omega 
\left\{-v\partial_t\varphi+\nabla v\cdot\nabla\varphi-cv\varphi\right\}\omega\,dzdt=0
$$
for all $\varphi\in C_0^\infty(\Omega\times I)$. 
The following results hold (see \cite{CS} and also \cite{I01}). 
\begin{proposition}
\label{Proposition:2.1}
Assume that $\omega$ is an $A_2$~weight on $B(0,1)$. 
Let $v$ be a solution of \eqref{eq:2.2} on $B(0,1)\times(0,1)$. 
Then there exists a constant $\gamma_1$ such that 
$$
\|v\|_{L^\infty(B(0,1/2)\times(1/2,1))}
\le\left(\frac{\gamma_1}{\omega(B(0,1))}\int_0^1\int_{B(0,1)} v^2\omega\,dzdt\right)^{\frac{1}{2}}. 
$$
Here $\gamma_1$ depends only on $N$, $\omega(B(0,1))$ and $\|c\|_{L^\infty(0,1:L^\infty(B(0,1)))}$. 
\end{proposition}
\begin{proposition}
\label{Proposition:2.2}
Assume that $\omega$ is an $A_2$~weight on $B(0,1)$. 
Let $v$ be a nonnegative solution of \eqref{eq:2.2} on $B(0,1)\times(-1,1)$. 
Then there exists a constant $\gamma_2$ such that 
$$
\sup_{Q_-} v\le \gamma_2\inf_{Q^+}v,
$$
where 
$$
Q_-:=B\left(0,\frac{1}{2}\right)\times\left(-\frac{3}{4},-\frac{1}{4}\right),
\qquad
Q_+:=B\left(0,\frac{1}{2}\right)\times\left(\frac{1}{4},\frac{3}{4}\right).
$$
Here $\gamma_2$ depends only on $N$, $\omega(B(0,1))$ and $\|c\|_{L^\infty(0,1:L^\infty(B(0,1)))}$. 
\end{proposition}
By Proposition~\ref{Proposition:2.2} we have:
\begin{lemma}
\label{Lemma:2.1} 
Let $R>0$ and $w$  an $A_2$~weight on $B(0,R)$. 
Let $v$ be a nonnegative solution of \eqref{eq:2.2} on $B(0,R)\times(0,T)$, 
where $0<T<\infty$.  
Then there exists a positive constant $C$ such that
\begin{equation}
\label{eq:2.3}
v(x_1,t_1)\le Cv(x_2,t_2)\exp\left(C\frac{|x_1-x_2|^2}{t_2-t_1}+\frac{t_2}{t_1}\right)
\end{equation}
for all $x_1$, $x_2\in B(0,R/2)$ and $0<t_1\le t_2\le T$. 
Here $C$ depends on $\omega(B(0,R))$ and $\|c\|_{L^\infty(-1,1:L^\infty(B(0,1)))}$. 
\end{lemma}
{\bf Proof.}
Let $x\in B(0,R/2)$ and $0<t<T$. 
Assume that 
$$
Q:=B(x,r)\times(t-r^2,t+r^2)\subset B(0,R)\times(0,T)
$$
for some $r>0$. 
Set 
$$
\tilde{v}(z,s):=v(x+rz,t+r^2s),
\quad
\tilde{\omega}(z):=\omega(x+rz),
\quad
\tilde{c}(z,s):=r^2c(x+rz,t+r^2s), 
$$
for $z\in B(0,1)$ and $s\in(-1,1)$. 
Then $\tilde{v}$ satisfies 
$$
\partial_s \tilde{v}
=\frac{1}{\tilde{\omega}}\mbox{div}_z\,(\tilde{\omega}(z)\nabla_z \tilde{v})+\tilde{c}\tilde{v}
\quad\mbox{in}\quad B(0,1)\times(-1,1). 
$$
Since $\tilde{\omega}(B(0,1))=\omega(B(0,r))\le\omega(B(0,R))$, 
by Proposition~\ref{Proposition:2.2} we can find a positive constant $c$, independent of $x$, $t$ and $r$, 
such that 
$$
\sup_{Q_-} \tilde{v}\le c\inf_{Q_+} \tilde{v},
$$
where $Q_+$ and $Q_-$ are as in Proposition~\ref{Proposition:2.2}. 
This implies that 
$$
\sup_{Q_-(x,t;r)} v\le c\inf_{Q_+(x,t;r)} v,
$$
where 
\begin{equation*}
\begin{split}
Q_-(x,t;r):= & B\left(x,\frac{r}{2}\right)\times\left(t-\frac{3}{4}r^2,t-\frac{1}{4}r^2\right),\\
Q_+(x,t;r):= & B\left(x,\frac{r}{2}\right)\times\left(t+\frac{1}{4}r^2,t+\frac{3}{4}r^2\right).
\end{split}
\end{equation*}
Then, similarly to \cite[Theorem~E]{A} and \cite[Theorem~2]{Moser}, 
we obtain \eqref{eq:2.3}. (See also \cite{IM1} and \cite{IM2}.)  
Thus Lemma~\ref{Lemma:2.1} follows. 
$\Box$
\section{Behavior of the harmonic function}\label{section:3}
In this section we study the behavior of positive harmonic functions 
for nonnegative Schr\"odinger operators and prove Theorem~\ref{Theorem:1.1}. 
In what follows, for $\lambda\in[\lambda_*,\infty)$, 
set $u^\pm_\lambda(r) :=r^{A^\pm(\lambda)}$ if $\lambda>\lambda_*$ and 
$$
u^+_\lambda(r):=r^{-\frac{N-2}{2}},
\qquad
u^-_\lambda(r):=r^{-\frac{N-2}{2}}|\log r|
$$
if $\lambda = \lambda_*$. 
Furthermore, we put $V_\lambda(r):=V(r)-\lambda r^{-2}$. 
\vspace{3pt}

We first study the behavior of solutions of (O) at $r=0$ and $r=\infty$. 
\begin{lemma}
\label{Lemma:3.1}
Let $V\in C((0,\infty))$. Assume that 
\begin{equation}
\label{eq:3.1}
\lim_{r\to 0}r^{2-\theta}\left|V_{\lambda_1}(r)\right|=0
\end{equation}
for some $\lambda_1\in[\lambda_*,\infty)$ and $\theta>0$. 
Then there exist solutions $U_*^\pm$ of {\rm (O)} such that
\begin{equation}
\label{eq:3.2}
\begin{split}
U_*^\pm(r) & =u^\pm_{\lambda_1}(r)+O(r^{\theta'}u_{\lambda_1}^\pm(r)),\\
(U_*^\pm)'(r) & =(u^\pm_{\lambda_1})'(r)+O(r^{-1+\theta'}u_{\lambda_1}^\pm(r)),
\end{split}
\end{equation}
as $r\to 0$, for some $\theta'\in(0,\theta]$. 
Furthermore, for any solution $w$ of {\rm (O)}, 
there exist constants $C_1$ and $C_2$ such that 
\begin{equation}
\label{eq:3.3}
w(r)=C_1U_*^+(r)+C_2U_*^-(r),\qquad r>0. 
\end{equation}
\end{lemma}
{\bf Proof.}
The proof is similar to \cite[Section~3]{IK00} but we give  details for the sake of completeness. 
We write $u^\pm=u^\pm_{\lambda_1}$ for simplicity. 

We first construct the solution $U^+_*$ of (O). 
Set $U^+_1(r):=u^+(r)$ and define 
$U^+_n$ $(n=2,3,\dots)$ inductively by 
\begin{equation}
\label{eq:3.4}
 U^+_{n+1}(r):=u^+(r)(1+F_n(r)),
\end{equation}
where 
$$
F_n(r):=\int_0^r s^{1-N}[u^+(s)]^{-2}\left(\int_0^s \tau^{N-1}u^+(\tau)V_{\lambda_1}(\tau)U^+_n(\tau)\,d\tau\right)\,ds. 
$$
Let $0<R<1$ and assume that 
\begin{equation}
\label{eq:3.5}
|U^+_n(r)|\le 2u^+(r)\quad\mbox{in}\quad (0,R]
\end{equation}
for some $n\in\{1,2,\dots\}$. 
Then it follows from \eqref{eq:1.3}, \eqref{eq:3.1} and \eqref{eq:3.5} that 
\begin{equation}
\label{eq:3.6}
\begin{split}
|F_n'(r)| & \le Cr^{1-N}[u^+(r)]^{-2}\int_0^r \tau^{N-1}\tau^{-2+\theta}[u^+(\tau)]^2\,d\tau\\
 & =Cr^{-1-\sqrt{D_{\lambda_1}}}\int_0^r \tau^{-1+\theta+\sqrt{D_{\lambda_1}}}\,d\tau\le Cr^{-1+\theta}
\end{split}
\end{equation}
for $r\in(0,R]$. 
Taking a sufficiently small $R>0$ if necessary, 
by \eqref{eq:3.4} and \eqref{eq:3.6} 
we have 
\begin{equation}
\label{eq:3.7}
|U^+_{n+1}(r)-u^+(r)|\le Cr^\theta u^+(r)\le u^+(r) 
\end{equation}
for $r\in(0,R]$. 
This implies that \eqref{eq:3.5} holds for $n=1,2,\dots$.  
Furthermore, we see that 
\eqref{eq:3.7} holds for $n=1,2,\dots$. 
Applying the successive approximation arguments on the existence of solutions 
to ordinary differential equations (see e.g., \cite[Chapter~1]{CoLevi}), 
we can find a function $U^+_*\in C((0,R])$ such that 
\begin{equation}
\label{eq:3.8}
|U^+_*(r)-u^+(r)|\le Cr^\theta u^+(r),
\qquad
U^+_*(r)=u^+(r)(1+F(r)),
\end{equation}
for $r\in(0,R]$, 
where 
$$
F(r):=\int_0^r s^{1-N}[u^+(s)]^{-2}\left(\int_0^s \tau^{N-1}u^+(\tau)V_{\lambda_1}(\tau)U^+_*(\tau)\,d\tau\right)\,ds.
$$
Similarly to \eqref{eq:3.6}, 
it follows that $|F'(r)|\le Cr^{-1+\theta}$ on $(0,R]$, which implies that 
\begin{equation}
\label{eq:3.9}
(U^+_*)'(r)-(u^+)'(r)=(u^+)'(r)F(r)+u^+(r)F'(r)
=O(r^{-1+\theta}u^+(r))
\end{equation}
as $r\to 0$. Furthermore, since
$$
(U^+_*)''+\frac{N-1}{r}(U^+_*)'-\frac{\lambda_1}{r^2}U^+_*
=V_{\lambda_1}U^+_*
\quad\mbox{in}\quad(0,R],
$$
$U^+_*$ satisfies (O) on $(0,R]$. 
By \eqref{eq:3.8} and \eqref{eq:3.9}, 
extending $U^+_*$ to the solution of (O) on $(0,\infty)$,  
we obtain the desired solution $U^+_*$ of (O). 

Next we construct the solution $U^-_*$ in the case $\lambda_1=\lambda_*$. 
We set $U^-_1(r)=u^-(r)$ and define 
$U^-_n$ $(n=2,3,\dots)$ inductively by 
$$
U^-_{n+1}(r):=u^-(r)+u^+(r)\tilde{F}_n(r),
$$
where 
$$
\tilde{F}_n(r):=\int_0^r s^{1-N}[u^+(s)]^{-2}\left(\int_0^s \tau^{N-1}u^+(\tau)V_{\lambda_1}(\tau)U^-_n(\tau)\,d\tau\right)\,ds. 
$$
Let $0<R<1$ and assume that 
\begin{equation}
\label{eq:3.10}
|U^-_n(r)|\le 2u^-(r)\quad\mbox{in}\quad (0,R]
\end{equation}
for some $n\in\{1,2,\dots\}$. 
Similarly to \eqref{eq:3.6}, 
by \eqref{eq:1.2} and \eqref{eq:3.10} we have 
\begin{equation}
\label{eq:3.11}
\begin{split}
|\tilde{F}_n'(r)| & \le Cr^{1-N}[u^+(r)]^{-2}\int_0^r \tau^{N-1}\tau^{-2+\theta}u^+(\tau)u^-(\tau)\,d\tau\\
 & =Cr^{-1}\int_0^r \tau^{-1+\theta}|\log\tau|\,d\tau\le Cr^{-1+\theta}|\log r|
\end{split}
\end{equation}
for $r\in(0,R]$. 
This implies that 
$$
|\tilde{F}(r)|\le Cr^\theta|\log r|,
\qquad
u^+(r)|\tilde{F}(r)|\le Cr^{-\frac{N-2}{2}+\theta}|\log r|
=Cr^\theta u^-(r)
$$
for $r\in(0,R]$. 
Then, by a similar argument as in the construction of $U^+_*$ 
we can find the desired solution~$U^-_*$ in the case $\lambda=\lambda_*$. 

Next we construct the solution $U^-_*$ in the case $\lambda_1>\lambda_*$. 
Let $\delta$ be a sufficiently small positive constant. 
We set 
$U^-_1(r):=u^-(r)$ and define 
$U^-_n$ $(n=2,3,\dots)$ inductively by 
$$
U^-_{n+1}(r):=u^-(r)(1+G_n(r)),
$$
where
$$
G_n(r):=
\int_0^r s^{1-N}[u^-(s)]^{-2}\left(\int^\delta_s \tau^{N-1}u^-(\tau)V_{\lambda_1}(\tau)U^-_n(\tau)\,d\tau\right)\,ds.
$$
Similarly to \eqref{eq:3.5}, we assume 
\begin{equation}
\label{eq:3.12}
|U^-_n(r)|\le 2u^-(r)\quad\mbox{in}\quad (0,\delta]
\end{equation}
for some $n\in\{1,2,\dots\}$. 
Since we can assume, without loss of generality, that $\theta<\sqrt{D_{\lambda_1}}$,  
by \eqref{eq:1.2} and \eqref{eq:3.12} we have 
\begin{equation*}
\begin{split}
|G'_n(r)| 
 & \le Cr^{1-N}[u^-(r)]^{-2}\int^\delta_r \tau^{N-1}\tau^{-2+\theta}[u^-(\tau)]^2\,d\tau\\
 & =Cr^{-1+\sqrt{D_{\lambda_1}}}\int^\delta_r \tau^{-1+\theta-\sqrt{D_{\lambda_1}}}\,d\tau
 \le Cr^{-1+\theta}
\end{split}
\end{equation*}
for $r\in(0,\delta]$. 
Then, taking a sufficiently small $\delta>0$ if necessary, 
we obtain 
$$
|U^-_{n+1}(r)-u^-(r)|\le Cr^\theta u^-(r)\le u^-(r)
$$
for $r\in(0,\delta]$. 
Repeating the above argument, 
we can find the desired solution~$U^-_*$ in the case $\lambda>\lambda_*$. 
Therefore, we obtain the desired solutions $U^\pm_*$ of (O). 
Furthermore, 
since $U^\pm_*$ are linearly independent, 
we see \eqref{eq:3.3}. 
Thus Lemma~\ref{Lemma:3.1} follows.
$\Box$\vspace{3pt}
\begin{lemma}
\label{Lemma:3.2}
Let $V\in C((0,\infty))$. Assume that 
\begin{equation}
\label{eq:3.13}
\lim_{r\to \infty}r^{2+\theta}\left|V_{\lambda_2}(r)\right|=0
\end{equation}
for some $\lambda_2\in[\lambda_*,\infty)$ and $\theta>0$. 
Then there exist solutions $U^\pm_{**}$ of {\rm (O)} such that 
\begin{equation}
\label{eq:3.14}
U_{**}^\pm(r)=u^{\pm}_{\lambda_2}(r)+O(r^{-\theta'}u^{\pm}_{\lambda_2}(r))
\end{equation}
as $r\to\infty$, for some $\theta'\in(0,\theta]$. 
Furthermore, for any solution $w$ of {\rm (O)}, 
there exist constants $C_1$ and $C_2$ such that 
\begin{equation}
\label{eq:3.15}
w(r)=C_1U_{**}^+(r)+C_2U_{**}^-(r),\qquad r>0. 
\end{equation}
\end{lemma}
{\bf Proof.}
Let $w$ be a solution of (O) on $(0,\infty)$. 
Set $\hat{w}(s):=s^{-N+2}w(s^{-1})$, 
which is the Kelvin transformation of $w$. 
Then $\hat{w}$ satisfies 
\begin{equation}
\label{eq:3.16}
\hat{w}''+\frac{N-1}{s}\hat{w}'-\hat{V}(s)\hat{w}=0\quad\mbox{in}\quad(0,\infty),
\end{equation}
where $\hat{V}(s):=s^{-4}V(s^{-1})$. 
It follows from \eqref{eq:3.13} that 
$$
\lim_{s\to 0}s^{-\theta}|s^2\hat{V}(s)-\lambda_2|
=\lim_{r\to\infty}r^{2+\theta}|V_{\lambda_2}(r)|=0.
$$
Therefore, by Lemma~\ref{Lemma:3.1} 
we can find solutions $W^\pm(r)$ of \eqref{eq:3.16} such that 
\begin{equation}
\label{eq:3.17}
W^\pm(s)=u^{\pm}_{\lambda_2}(s)+O\left(s^{\theta'}u^{\pm}_{\lambda_2}(s)\right)\quad\mbox{as}\quad s\to 0,
\end{equation}
for some $\theta'\in(0,\theta]$. 
Set $U_{**}^\pm(r):=r^{-N+2}W^\mp(r^{-1})$. 
Then $U_{**}^\pm(r)$ are solutions of (O) on $(0,\infty)$. 
Furthermore, it follows that
$$
-A^\pm(\lambda_2)-N+2
=\frac{-(N-2)\mp\sqrt{(N-2)^2+4\lambda_2}}{2}
=A^\mp(\lambda_2),
$$
which together with \eqref{eq:3.17} implies \eqref{eq:3.14}.  
Furthermore, since $U^\pm_{**}$ are linearly independent, 
we obtain \eqref{eq:3.15}. 
Thus Lemma~\ref{Lemma:3.2} follows. 
$\Box$
\vspace{3pt}

In what follows, we set 
$$
U(r):=U^+_*(r). 
$$
Next we show the positivity of $U$ under the assumption that $H$ is nonnegative. 
\begin{lemma}
\label{Lemma:3.3}
Let $V\in C((0,\infty))$. Assume \eqref{eq:1.2} and that $H$ is nonnegative. 
Then $U(r)>0$ on $(0,\infty)$. 
\end{lemma}
{\bf Proof.} 
We consider the case $-\lambda_*\le\lambda_1\le 0$. 
For $n=1,2,\dots$, set 
$$
V_n(r):=\max\{-n,V(r)\},
\qquad
H_n:=-\Delta+V_n. 
$$
Since $V_n\in L^\infty(0,\infty)$ and $H_n$ is nonnegative, 
by (ii) of Theorem 3.1 in \cite{M01} there exists a radially symmetric and bounded function $u_n=u_n(|x|)\in C^2({\bf R}^N)$ such that 
\begin{equation}
\label{eq:3.18}
-\Delta u_n+V_nu_n=0\quad\mbox{in}\quad{\bf R}^N,
\qquad
u_n>0\quad\mbox{in}\quad{\bf R}^N. 
\end{equation}
In particular, it follows from the regularity of $u_n$ that 
\begin{equation}
\label{eq:3.19}
u_n'(0)=0.
\end{equation}
By \eqref{eq:3.2} we can find $R>0$ such that $U(r)>0$ on $(0,R]$. 
Set
$$
U_n(r):=U(R)\frac{u_n(r)}{u_n(R)}. 
$$
Since $U_n$ satisfies \eqref{eq:3.18}, 
we have 
$$
-\frac{1}{r^{N-1}}(r^{N-1}U_n')'+V_nU_n=0\quad\mbox{in}\quad(0,\infty), 
$$
which implies that 
\begin{equation}
\label{eq:3.20}
\begin{split}
0 & =\int_{r'}^r[-(s^{N-1}U_n')'U+s^{N-1}V_nU_nU]\,ds\\
 & =\left[-s^{N-1}U_n'U\right]_{s=r'}^{s=r}
+\int_{r'}^r[s^{N-1}U_n'U'+s^{N-1}V_nU_nU]\,ds
\end{split}
\end{equation}
for $0<r'<r$. 
Similarly, since $U$ is a solution of (O), we have 
\begin{equation}
\label{eq:3.21}
\begin{split}
0 & =\int_{r'}^r[-(s^{N-1}U')'U_n+s^{N-1}V(s)UU_n]\,ds\\
 & =\left[-s^{N-1}U'U_n\right]_{s=r'}^{s=r}
+\int_{r'}^r[s^{N-1}U'U_n'+s^{N-1}V(s)UU_n]\,ds
\end{split}
\end{equation}
for $0<r'<r$. 
Since $U(r)>0$, $U_n(r)>0$ and $V(r)\le V_n(r)$ on $(0,R]$, 
we deduce from \eqref{eq:3.20} and \eqref{eq:3.21} that
\begin{equation}\label{eq:3.21_1}
\begin{split}
 & r^{N-1}[U'(r)U_n(r)-U(r)U_n'(r)]
-(r')^{N-1}[U'(r')U_n(r')-U(r')U_n'(r')]\\
 & =\int_{r'}^rs^{N-1}[V(s)-V_n(s)]UU_n\,ds\le 0
\end{split}
\end{equation}
for $0<r'<r\le R$. 
On the other hand, 
it follows from \eqref{eq:3.2} that 
$$
\lim_{r'\to 0}(r')^{N-1}[U'(r')U_n(r')-U(r')U_n'(r')]=0. 
$$
Taking $r'\to 0$ in \eqref{eq:3.21_1} 
together with \eqref{eq:3.19} implies that 
$$
0\ge r^{N-1}[U'(r)U_n(r)-U(r)U_n'(r)]
=r^{N-1}U_n(r)^2\left(\frac{U(r)}{U_n(r)}\right)',
\qquad 0<r\le R.
$$
We deduce from $U_n(R)=U(R)$ that 
$$
\frac{U(r)}{U_n(r)}\ge\frac{U(R)}{U_n(R)}=1, \qquad 0<r\le R.
$$
Therefore we obtain 
\begin{equation}
\label{eq:3.22}
0<U_n(r)\le U(r),\qquad 0<r\le R.
\end{equation}

On the other hand, 
since $V\in C((0,\infty))$ and $U_n(R)=1$, 
by the Harnack inequality and regularity theorems for elliptic equations 
in a similar way to the Perron method, 
we can find a function $\tilde{U}\in C^2((0,\infty))$ such that 
$$
\lim_{n\to\infty}\|U_n-\tilde{U}\|_{C^2(I)}=0
$$
for any compact set $I$ in $(0,\infty)$. 
Then $\tilde{U}$ is a solution of (O) on $(0,\infty)$. 
Furthermore, by \eqref{eq:3.22} we see that 
\begin{equation}
\label{eq:3.23}
\tilde{U}(r)\ge 0\quad\mbox{in}\quad(0,\infty),
\qquad
\tilde{U}(r)\le U(r)\quad\mbox{in}\quad(0,R],
\qquad
\tilde{U}(R)=U(R)>0. 
\end{equation}
Using the Harnack inequality again, we obtain
\begin{equation}
\label{eq:3.24}
\tilde{U}(r)>0 \quad\mbox{in}\quad(0,\infty).
\end{equation}
Furthermore, 
by Lemma~\ref{Lemma:3.1} there exist constants $C_1$ and $C_2$ such that 
\begin{equation}
\label{eq:3.25}
\tilde{U}(r)=C_1U_*^+(r)+C_2U^-_*(r)
\quad\mbox{in}\quad(0,\infty). 
\end{equation}
Since $A^-(\lambda_1)<A^+(\lambda_1)$, 
by \eqref{eq:3.23} and \eqref{eq:3.25} 
we see that $C_2=0$ and $C_1=1$, 
that is $\tilde{U}(r)=U^+_*(r)=U(r)$ on $(0,\infty)$. 
Therefore we deduce from \eqref{eq:3.24} that $U(r)>0$ on $(0,\infty)$. 

It remains to consider the case $\lambda_1>0$. 
Let $k\in\{1,2,\dots\}$ be such that $\lambda_1<\omega_k:=k(N+k-2)$. 
For any $\phi\in C^\infty({\bf R}^{N+2k}\setminus\{0\})$ and $\omega\in{\bf S}^{N+2k-1}$, 
set 
$$
\phi_\omega(r):=\phi(r\omega),
\qquad
\psi_\omega(r):=r^k\phi_\omega(r).
$$
Since $H$ is nonnegative, 
we have 
\begin{equation*}
\begin{split}
0 & \le\frac{1}{|{\bf S}^{N-1}|}\int_{{\bf R}^N}\left[|\nabla\psi_\omega|^2+V\psi_\omega^2\right]\,dx
 =\int_0^\infty\left[|\psi_\omega'|^2+V\psi_\omega^2\right]r^{N-1}\,dr\\
 & =\int_0^\infty
\left[r^{2k}|\phi_\omega'|^2+2kr^{2k-1}\phi_\omega\phi_\omega'+k^2r^{2k-2}\phi_\omega^2+r^{2k}V\phi_\omega^2\right]r^{N-1}\,dr\\
 & =\int_0^\infty
\left[|\phi_\omega'|^2+k^2r^{-2}\phi_\omega^2+V\phi_\omega^2\right]r^{N+2k-1}\,dr
+\int_0^\infty kr^{N+2k-2}[(\phi_\omega)^2]'\,dr\\
 & =\int_0^\infty
 \left[|\phi_\omega'|^2+\tilde{V}\phi_\omega^2\right]r^{N+2k-1}\,dr
 \le|{\bf S}^{N+2k-1}|^{-1}\int_{{\bf R}^{N+2k}}
 \left[|\nabla\phi|^2+\tilde{V}\phi^2\right]\,dx,
\end{split}
\end{equation*}
where $\tilde{V}(r):=V(r)-\omega_kr^{-2}$. 
This means that 
$\tilde{H}:=-\Delta_{N+2k}+\tilde{V}$
is nonnegative operator on ${\bf R}^{N+2k}$. 
Furthermore, \eqref{eq:1.2} holds with $\lambda_1$ and $\lambda_2$ replaced by 
$$
\tilde{\lambda}_1:=\lambda_1-\omega_k
>-\frac{(N+2k-2)^2}{4}\quad\mbox{and}\quad
\tilde{\lambda}_2:=\lambda_2-\omega_k
>-\frac{(N+2k-2)^2}{4},
$$
respectively. 
Therefore, by Lemma~\ref{Lemma:3.3} in the case $\lambda_*<\lambda_1\le 0$ 
we can find a solution $u=u(r)$ of 
\begin{equation*}
\begin{split}
 & u''+\frac{N+2k-1}{r}u'-\tilde{V}(r)=0\quad\mbox{in}\quad(0,\infty),\\
 & u(r)>0\quad\mbox{in}\quad(0,\infty),
\quad
u(r)=r^a+O(r^{a+\theta})\quad\mbox{as}\quad r\to 0, 
\end{split}
\end{equation*}
where 
$$
a:=\frac{-(N+2k-2)+\sqrt{(N+2k-2)^2+4(\lambda_1-\omega_k)}}{2}
=-k+A^+(\lambda_1). 
$$
Then $\tilde{U}(r):=r^ku(r)$ is a solution of (O) and it satisfies
$$
\tilde{U}(r)>0\quad\mbox{in}\quad(0,\infty),
\qquad
\tilde{U}(r)\thicksim r^{A^+(\lambda_1)}
\quad\mbox{as}\quad r\to 0.
$$
It follows from Lemma~\ref{Lemma:3.1} that 
$U(r)=\tilde{U}(r)>0$ on $(0,\infty)$. 
Thus Lemma~\ref{Lemma:3.3} follows.
$\Box$
\vspace{5pt}

Next we study the asymptotic behavior of $U(r)$ as $r\to\infty$. 
\begin{lemma}
\label{Lemma:3.4}
Let $V\in C((0,\infty))$. Assume \eqref{eq:1.2} and that 
$-\Delta+V(|x|)-W(|x|)$ is nonnegative for some $W\in C_0((0,\infty))$ with 
$$
W\ge 0,\qquad W\not\equiv 0\quad\mbox{in}\quad (0,\infty).
$$ 
Then there exists a positive constant $c$ such that 
$$
U(r)\thicksim cr^{A^+(\lambda_2)}\,\,\,\mbox{if $\lambda_2>\lambda_*$},
\qquad
U(r)\thicksim cr^{-\frac{N-2}{2}}\log r\,\,\,\mbox{if $\lambda_2=\lambda_*$},
$$
as $r\to\infty$.
\end{lemma}
{\bf Proof.}
Since $-\Delta+V - W$ is nonnegative, 
by Lemmas~\ref{Lemma:3.1} and \ref{Lemma:3.3}
we can find a function $U_W\in C^2((0,\infty))$ satisfying
\begin{equation}
\label{eq:3.26}
\begin{split}
 & U_W''+\frac{N-1}{r}U_W'-(V(r)-W(r))U_W=0\quad\mbox{in}\quad(0,\infty),\\
 & U_W>0\quad\mbox{in}\quad(0,\infty),
\qquad
U_W(r)=r^{A^+(\lambda_1)}(1+o(1))\quad\mbox{as}\quad r\to 0.
\end{split}
\end{equation}
On the other hand, $U$ satisfies
$$
U''+\frac{N-1}{r}U'-(V(r)-W(r))U=W(r)U\quad\mbox{in}\quad(0,\infty). 
$$
Define 
$$
\tilde{U}(r)=U_W(r)(1+F_W(r)),
$$
where 
$$
F_W(r):=\int_0^r s^{1-N}[U_W(s)]^{-2}\left(\int_0^s \tau^{N-1}U_W(\tau)W(\tau)U(\tau)\,d\tau\right)\,ds\ge 0.
$$
Since $W$ has a compact support, 
$F_W(r)=0$ for all sufficiently small $r>0$. 
Furthermore, by \eqref{eq:3.26} we have 
$$
\tilde{U}''+\frac{N-1}{r}\tilde{U}'-(V(r)-W(r))\tilde{U}
=WU\quad\mbox{in}\quad(0,\infty). 
$$
Then $\hat{U}:=U-\tilde{U}$ satisfies 
\begin{equation*}
\begin{split}
 & \hat{U}''+\frac{N-1}{r}\hat{U}'-(V(r)-W(r))\hat{U}=0\quad\mbox{in}\quad(0,\infty),\\
 & \hat{U}(r)=U(r)-U_W(r)=o\left(r^{A^+(\lambda_1)}\right)
\quad\mbox{as}\quad r\to 0.
\end{split}
\end{equation*}
This  together with Lemma~\ref{Lemma:3.1} imply that 
$\hat{U}=0$ in $(0,\infty)$, that is, 
\begin{equation}
\label{eq:3.27}
U(r)=U_W(r)(1+F_W(r))\quad\mbox{in}\quad(0,\infty).
\end{equation}
On the other hand, by Lemmas~\ref{Lemma:3.2} and \ref{Lemma:3.3} 
we see that, either 
$$
{\rm (a)}\quad U_W(r)\thicksim c_1u^-_{\lambda_2}(r)
\qquad\mbox{or}
\qquad
{\rm (b)}\quad U_W(r)\thicksim c_2u^+_{\lambda_2}(r)
$$
as $r\to\infty$, where $c_1$ and $c_2$ are positive constants. 

Consider the case $\lambda_2>\lambda_1$. 
Assume that (a) holds. 
Since $W$ has a compact support and $U>0$ on $(0,\infty)$, 
we take a sufficiently large constant $R>0$ so that 
\begin{equation}
\label{eq:3.28}
\begin{split}
 & U_W(r)F_W(r)\ge CU_W(r)\int_R^r s^{1-N}[U_W(s)]^{-2}\,ds\\
 & \qquad
\ge Cr^{\frac{-(N-2)-\sqrt{D(\lambda_2)}}{2}}
\int_1^r s^{-1+\sqrt{D(\lambda_2)}}\,ds\ge Cr^{\frac{-(N-2)+\sqrt{D(\lambda_2)}}{2}}
=Cr^{A^+(\lambda_2)}
\end{split}
\end{equation}
for all sufficiently large $r$.  
On the other hand, if (b) holds, then 
it follows from \eqref{eq:3.27} that 
\begin{equation}
\label{eq:3.29}
U(r)\ge U_W(r)\ge Cr^{A^+(\lambda_2)}
\end{equation}
for all sufficiently large $r$. 
In both cases of (a) and (b), 
$U(r)\ge Cr^{A^+(\lambda_2)}$ for all sufficiently large $r$. 
Then Lemma~\ref{Lemma:3.4} in the case~$\lambda_2>\lambda_*$ follows from Lemma~\ref{Lemma:3.2}. 

Consider the case $\lambda_2=\lambda_*$. 
If (a) holds, then, similarly to \eqref{eq:3.29}, we have 
$U(r)\ge U_W(r)\ge Cr^{-\frac{N-2}{2}}\log r$ 
for all sufficiently large $r$. 
If (b) holds, then, similarly to \eqref{eq:3.28}, 
we have 
\begin{equation*}
U_W(r)F_W(r)\ge CU_W(r)\int_R^r s^{1-N}[U_W(s)]^{-2}\,ds
\ge Cr^{-\frac{N-2}{2}}
\int_1^r s^{-1}\,ds\ge Cr^{-\frac{N-2}{2}}\log r
\end{equation*}
for all sufficiently large $r$. 
In both cases of (a) and (b), 
$U(r)\ge Cr^{-\frac{N-2}{2}}\log r$ for all sufficiently large $r$. 
Then Lemma~\ref{Lemma:3.4} in the case~$\lambda_2=\lambda_*$ follows from Lemma~\ref{Lemma:3.2}. 
Thus the proof of Lemma~\ref{Lemma:3.4} is complete. 
$\Box$\vspace{3pt}
\newline
Next we employ the arguments in \cite[Lemma~6]{DS} 
and prove the following lemma. 
\begin{lemma}
\label{Lemma:3.5}
Let $V\in C((0,\infty))$. Assume \eqref{eq:1.2} and that $H$ is nonnegative. 
If there exists a positive constant $c$ such that 
\begin{equation}
\label{eq:3.30}
U(r)\ge cr^{A^+(\lambda_2)}\quad\mbox{if}\quad\lambda_2>\lambda_*,
\qquad
U(r)\ge cr^{-\frac{N-2}{2}}\log r\quad\mbox{if}\quad\lambda_2=\lambda_*,
\end{equation}
for all sufficiently large $r>0$, then $H:=-\Delta+V$ is subcritical. 
\end{lemma}
{\bf Proof.}
Let $\phi\in C_0^\infty({\bf R}^N\setminus\{0\})$. Set 
$$
\tilde{\phi}(x):=\phi(x)/U(|x|). 
$$
Then we have
\begin{equation*}
\begin{split}
|\nabla\tilde{\phi}|^2
 & =\left|\frac{U\nabla\phi-\phi\nabla U}{U^2}\right|^2
=\frac{U^2|\nabla\phi|^2-2U\phi\nabla\phi\nabla U+\phi^2|\nabla U|^2}{U^4}\\
 & =\frac{U^2|\nabla\phi|^2-U\nabla|\phi|^2\nabla U+\phi^2|\nabla U|^2}{U^4}.
\end{split}
\end{equation*}
This implies that 
\begin{equation}
\label{eq:3.31}
\begin{split}
\int_{{\bf R}^N}
|\nabla\tilde{\phi}|^2U^2\,dx
 & =\int_{{\bf R}^N}\left[|\nabla\phi|^2+\phi^2\nabla\left(\frac{\nabla U}{U}\right)
 +\frac{\phi^2|\nabla U|^2}{U^2}\right]U^2\,dx\\
 & =\int_{{\bf R}^N}(|\nabla\phi|^2+V\phi^2)\,dx.
\end{split}
\end{equation}

Let $\tilde{V}$ satisfy \eqref{eq:1.2} and 
$\tilde{V}(r)>\lambda_*r^{-2}$ on $(0,\infty)$.
Let $W\in C_0^\infty({\bf R}^N)$ be such that $W\ge 0$ on ${\bf R}^N$.
By the Hardy inequality there exists $\epsilon>0$ such that 
\begin{equation}
\label{eq:3.32}
\int_{{\bf R}^N}
\left[|\nabla\psi|^2+(\tilde{V}-\epsilon W)\psi^2\right]\,dx\ge 0
\quad\mbox{for}\quad
\psi\in C^\infty({\bf R}^N\setminus\{0\}). 
\end{equation}
Furthermore, by Lemma~\ref{Lemma:3.4} 
we can find a positive function $\tilde{U}\in C^2((0,\infty))$ such that 

\begin{equation}
\label{eq:3.33}
\begin{split}
 & \tilde{U}''+\frac{N-1}{r}\tilde{U}'-\tilde{V}\tilde{U}=0\quad\mbox{in}\quad(0,\infty),
 \qquad\tilde{U}(r)\thicksim r^{A^+(\lambda_1)}\quad\mbox{as}\quad r\to 0,\\
 & \tilde{U}(r)\thicksim cr^{A^+(\lambda_2)}\quad\mbox{as}\quad r\to \infty\quad\mbox{if}\quad\lambda>\lambda_*,\\
 & \tilde{U}(r)\thicksim cr^{-\frac{N-2}{2}}\log r\quad\mbox{as}\quad r\to \infty\quad\mbox{if}\quad\lambda=\lambda_*,
\end{split}
\end{equation}
for some constant $c>0$. 
Since $U(r)=r^{A^+(\lambda_1)}(1+o(1))$ as $r\to 0$, 
it follows from \eqref{eq:3.30} and \eqref{eq:3.33} that 
$U(r)\ge C^{-1}\tilde{U}(r)$ on $(0,\infty)$. 
This together with \eqref{eq:3.32} implies that 
\begin{equation}
\label{eq:3.34}
\begin{split}
\int_{{\bf R}^N}
|\nabla\tilde{\phi}|^2U^2\,dx
 & \ge C\int_{{\bf R}^N}
|\nabla\tilde{\phi}|^2\tilde{U}^2\,dx\\
 & =C\int_{{\bf R}^N}(|\nabla\hat{\phi}|^2+\tilde{V}\hat{\phi}^2)\,dx
 \ge C\epsilon\int_{{\bf R}^N}W\hat{\phi}^2\,dx,
\end{split}
\end{equation}
as in the same way as \eqref{eq:3.31}, 
where
$$
\hat{\phi}(x)=\tilde{U}(|x|)\tilde{\phi}(x)=\frac{\tilde{U}(|x|)}{U(|x|)}\phi(x).
$$
Since $U(r)\thicksim\tilde{U}(r)$ as $r\to 0$ and $W$ has a compact support, 
we deduce from \eqref{eq:3.31} and \eqref{eq:3.34} that
$$
\int_{{\bf R}^N}(|\nabla\phi|^2+V\phi^2)\,dx
\ge C\epsilon\int_{{\bf R}^N}W\hat{\phi}^2\,dx
\ge C\epsilon\int_{{\bf R}^N}W\phi^2\,dx,
\quad \phi\in C_0^\infty({\bf R}^N\setminus\{0\}).
$$
This means that $H$ is subcritical. 
Thus Lemma~\ref{Lemma:3.5} follows. 
$\Box$
\vspace{5pt}
\newline
{\bf Proof of Theorem~\ref{Theorem:1.1}.}
Assertions~(1) and (2) follow from Lemma~\ref{Lemma:3.1}. 
Assertion~(3) follows from Lemmas~\ref{Lemma:3.2}, \ref{Lemma:3.3} and \ref{Lemma:3.5}. 

It remains to prove assertion~(4). 
Let $W\in C_0((0,\infty))$ be such that $W\ge 0$ and $W\not\equiv 0$ on $(0,\infty)$. 
Assume that $H$ is nonnegative. 
For any $\mu\in{\bf R}$, let
$$
H_\mu:=-\Delta+V+\mu W. 
$$
Define 
$$
I:=\{\mu\in{\bf R}\,:\,\mbox{$H_\mu$ is subcritical}\},
\qquad
\mu_*:=\inf_{\mu\in I}\mu. 
$$
It follows from Lemmas~\ref{Lemma:3.4} and \ref{Lemma:3.5} that
$I=(\mu_*,\infty)$ and $\mu_*\le 0$.  
Since $H_{\mu_*}$ is nonnegative, $H_{\mu_*}$ must be critical. 
Then assertion~(4) follows. 
Therefore the proof of Theorem~\ref{Theorem:1.1} is complete
by replacing  $W$ in this proof by $-W$. 
$\Box$
\section{Proof of Theorem~\ref{Theorem:1.2}}
Assume \eqref{eq:1.2}. 
Let $H:=-\Delta+V$ be nonnegative 
and $U$ the positive harmonic function 
given in Theorem~\ref{Theorem:1.1}. 
We define the unitary operator ${\mathcal U}$ by 
$$
{\mathcal U}\,:\, L^2({\bf R}^N,dx)\ni f\,\,\longmapsto\,\,U^{-1}f\in L^2({\bf R}^N,\omega(x)\,dx),
$$
where $\omega(x)=U(|x|)^2$. 
Then the operator $L := {\mathcal U} H {\mathcal U}^{-1}$ is given by 
$$
Lv:=-\frac{1}{\omega(x)}\mbox{div}\,(\omega(x)\nabla v).
$$ 
We denote by $p(x,y,t)$ and $G(x,y,t)$ the heat kernels of $H$ and $L$, respectively. 
Then 
\begin{equation}
\label{eq:4.1}
p(x,y,t) = U(|x|)U(|y|)G(x,y,t)
\end{equation}
for $x$, $y\in{\bf R}^N$ and $t>0$. 
In this section we study upper and lower bounds of $G=G(x,y,t)$
and then obtain Theorem \ref{Theorem:1.2}. 
\begin{lemma}
\label{Lemma:4.1}
Let $x$, $y\in{\bf R}^N$ and $t>0$. 
Assume that $\omega$ is an $A_2$ weight on $B(x,2\sqrt{t})\cup B(y,2\sqrt{t})$. 
Then there exists a constant $C$ such that 
\begin{equation}
\label{eq:4.2}
G(x,y,t)\le
\frac{C}{\sqrt{\omega(B(x,\sqrt{t}))}\sqrt{\omega(B(y,\sqrt{t}))}}\exp\left(-\frac{|x-y|^2}{Ct}\right)
\end{equation}
for $x$, $y\in{\bf R}^N$ and $t>0$, where  
$C = C_{x,y,t}$ depends on 
$[\omega](B(x, 2\sqrt{t}))$ and $[\omega](B(y, 2\sqrt{t}))$. 
In particular, $C$ is independent of 
$x$, $y$ and $t$ if $w$ is an $A_2$~weight on ${\bf R}^N$.
\end{lemma}
{\bf Proof.}
We obtain the upper bound of $G=G(x,y,t)$ by using the standard method as given, 
for example, \cite[Section~6]{S}. 
We give the proof for completeness of this paper. 

We fix  $x$, $y\in{\bf R}^N$ and $t>0$. 
Let $\lambda\in{\bf R}$ and  $\psi$ be a bounded smooth function  on ${\bf R}^N$ such that $|\nabla\psi|\le 1$ on ${\bf R}^N$. 
For $f_0\in L^2(B(y,\sqrt{t}),\omega\,dz)$, we set 
\begin{equation}
\label{eq:4.3}
 f(\xi,s):=\int_{B(y,\sqrt{t})}G(\xi,z,s)e^{-\lambda\psi(z)}f_0(z)\omega(z)\,dz,
 \quad
 F(\xi,s):=e^{\lambda\psi(\xi)}f(\xi,s).
\end{equation}
Since $f=f(\xi,s)$ satisfies
\begin{equation}
\label{eq:4.4}
\partial_s f=\frac{1}{\omega(\xi)}\mbox{div}_\xi\,(\omega(\xi)\nabla_\xi f)\quad\mbox{in}\quad{\bf R}^N\times(0,\infty), 
\end{equation}
we have 
\begin{equation*}
\begin{split}
 \frac{d}{ds}\int_{{\bf R}^N}F(\xi,s)^2\omega\,d\xi
 & =2\int_{{\bf R}^N}e^{2\lambda\psi(\xi)}f(\partial_s f)\omega\,d\xi\\
 & =-2\int_{{\bf R}^N}\left[2\lambda e^{2\lambda\psi}f\,\nabla\psi\cdot\nabla f
 +e^{2\lambda\psi}|\nabla f|^2\right]\,\omega\,d\xi\\
 & \le 2\lambda^2\int_{{\bf R}^N}e^{2\lambda\psi}f^2|\nabla\psi|^2\omega\,d\xi\
 \le 2\lambda^2\int_{{\bf R}^N}F(\xi,s)^2\omega\,d\xi,
\end{split}
\end{equation*}
which implies that 
\begin{equation}
\label{eq:4.5}
\int_{{\bf R}^N}F(\xi,s)^2\omega\,d\xi\le e^{2\lambda^2 s}\int_{{\bf R}^N}F(\xi,0)^2\omega\,d\xi
=e^{2\lambda^2 s}\int_{B(y,\sqrt{t})}f_0(\xi)^2\omega\,d\xi,
\quad s>0. 
\end{equation}

Let $0<\tau\le t$ and let $\tilde{f}(y,s):=f(x+\sqrt{\tau}y/2,3\tau/4+\tau s/4)$ for $(y,s)\in B(0,1)\times(0,1)$. 
Then $\tilde{f}$ satisfies 
$$
\partial_s\tilde{f}=\frac{1}{\tilde{\omega}}\mbox{div}\,(\tilde{\omega}\nabla\tilde{f})
\quad\mbox{in}\quad B(0,1)\times(0,1),
$$
where $\tilde{\omega}(y):=\omega(x+\sqrt{\tau}y/2)$. 
Then, by Proposition~\ref{Proposition:2.1} we obtain 
\begin{equation}
\label{eq:4.6}
\begin{split}
f(x,\tau)^2=\tilde{f}(0,0)^2
 & \le\frac{C}{\tilde{\omega}(B(0,1))}\int_0^1\int_{B(0,1)}\tilde{f}(y,s)^2\tilde{\omega}\,dyds\\
 & \le\frac{C}{\tau \omega(B(x,\sqrt{\tau}))}\int_{3\tau/4}^\tau\int_{B(x,\sqrt{\tau})}f(\xi,s)^2\omega\,d\xi ds. 
\end{split}
\end{equation}
The constant $C$ depends  on $[\omega](B(x, 2\sqrt{t}))$. 
Since $|\nabla\psi|\le 1$, 
by \eqref{eq:2.1}, \eqref{eq:4.3}, \eqref{eq:4.5} and \eqref{eq:4.6} we have 
\begin{equation}
\label{eq:4.7}
\begin{split}
e^{2\lambda\psi(x)}f(x,\tau)^2
 & \le\frac{C}{\tau \omega(B(x,\sqrt{t}))}
 \int_{3\tau/4}^\tau\int_{B(x,\sqrt{\tau})}e^{2\lambda(\psi(x)-\psi(\xi))}F(\xi,s)^2\omega\,d\xi ds\\
 & \le\frac{Ce^{2\lambda\sqrt{t}}}{\omega(B(x,\sqrt{t}))}\sup_{0<s<t}\int_{{\bf R}^N}F(\xi,s)^2\omega\,d\xi\\
 & \le\frac{Ce^{2\lambda\sqrt{t}}}{\omega(B(x,\sqrt{t}))}e^{2\lambda^2 t}
 \int_{B(y,\sqrt{t})}f_0(\xi)^2\omega\,d\xi
\end{split}
\end{equation}
for all $t/2\le\tau\le t$. 
Furthermore, by \eqref{eq:4.3} we obtain 
\begin{equation*}
\begin{split}
f(x,\tau) & =e^{-\lambda\psi(y)}\int_{B(y,\sqrt{t})}G(x,z,\tau)e^{-\lambda(\psi(z)-\psi(y))}f_0(z)\omega\,dz\\
 & \ge e^{-\lambda\sqrt{t}}e^{-\lambda\psi(y)}
 \int_{B(y,\sqrt{t})}G(x,z,\tau)f_0(z)\omega\,dz
\end{split}
\end{equation*}
for $\tau>0$. 
This implies that 
\begin{equation*}
\begin{split}
 & \left(\int_{B(y,\sqrt{t})}G(x,z,\tau)^2\omega(z)\,dz\right)^{1/2}\\
 & =\sup\,\biggr\{\int_{{\bf R}^N}G(x,z,\tau)f_0(z)\omega(z)\,dz\,:\,\\
 & \qquad\qquad\qquad
 f_0\in L^2(B(y,\sqrt{t}),\omega\,dz),\,\,\,
 \|f_0\|_{L^2(B(y,\sqrt{t}),\omega\,dz)}\le 1\biggr\}\\
 & \le e^{\lambda\sqrt{t}}e^{\lambda\psi(y)}\sup\left\{f(x,\tau)\,:\,f_0\in L^2(B(y,\sqrt{t}),\omega\,dz),\,\,\,
 \|f_0\|_{L^2(B(y,\sqrt{t}),\omega\,dz)}\le 1\right\},
\end{split}
\end{equation*}
which together with \eqref{eq:4.7} yields
\begin{equation}
\label{eq:4.8}
\begin{split}
 & \int_{B(y,\sqrt{t})}G(x,z,\tau)^2\omega(z)\,dz\\
 & \le e^{2\lambda\sqrt{t}}e^{-2\lambda(\psi(x)-\psi(y))}\\
 &\qquad\times
 \sup\left\{e^{2\lambda\psi(x)}f(x,\tau)^2\,:\,f_0\in L^2(B(y,\sqrt{t}),\omega\,dz),\,\,\,
 \|f_0\|_{L^2({\bf R}^N,\omega\,dz)}\le 1\right\}\\
 & \le e^{2\lambda\sqrt{t}}e^{-2\lambda(\psi(x)-\psi(y))}
 \frac{Ce^{2\lambda\sqrt{t}}}{\omega(B(x,\sqrt{t}))}e^{2\lambda^2 t}
\end{split}
\end{equation}
for all $t/2\le\tau\le t$. 

On the other hand, 
since $\tilde{g}(\xi,s):=G(x,\xi,s)$ is also a solution of \eqref{eq:4.4}, 
similarly to \eqref{eq:4.6}, 
we have 
$$
G(x,y,t)^2\le\frac{C}{t\omega(B(y,\sqrt{t}))}
 \int_{3t/4}^t\int_{B(y,\sqrt{t})}G(x,z,\tau)^2\omega(z)\,dzd\tau.
$$
Then we deduce from \eqref{eq:4.8} that 
$$
G(x,y,t)^2
\le\frac{C}{\omega(B(x,\sqrt{t}))\omega(B(y,\sqrt{t}))}
e^{4\lambda\sqrt{t}+2\lambda^2t-2\lambda(\psi(x)-\psi(y))}.
$$ 
We choose  $\lambda = \frac{\psi(x) - \psi(y)}{2t}$ and optimize over $\psi$ with  $|\nabla \psi | \le 1$. 
This gives \eqref{eq:4.2}, and the proof is complete. 
$\Box$\vspace{3pt}

If   $w$ is an $A_2$~weight  on ${\bf R}^N$, 
then we obtain upper estimate of Lemma~\ref{Lemma:4.1}. 
We mention that the proof actually gives  the estimate
$$
G(x,y,t)\le
\frac{C_\epsilon}{\sqrt{\omega(B(x,\sqrt{t}))}\sqrt{\omega(B(y,\sqrt{t}))}}\exp\left(-\frac{|x-y|^2}{(4+\epsilon)t}\right)
$$
for every $\epsilon > 0$ and all $x, y \in {\bf R}^N$ and $t > 0$.  Here $C_\epsilon$ is a positive constant depending on $\epsilon$. Therefore, by  \eqref{eq:4.1} we obtain  the following upper estimate  
\begin{equation}
\label{eq:4.9}
p(x,y,t)\le
\frac{C_\epsilon U(x)U(y)}{\sqrt{\omega(B(x,\sqrt{t}))}\sqrt{\omega(B(y,\sqrt{t}))}}\exp\left(-\frac{|x-y|^2}{(4+\epsilon)t}\right).
\end{equation}
This shows  the upper bound of Theorem \ref{Theorem:1.2}. 

Next,   we prove the lower bound. In the rest of this section we assume that $w$ is an $A_2$~weight on ${\bf R}^N$. The idea of proof is known and has been used in the context of Riemmannian manifolds, see, e.g., 
\cite{Cou}, \cite[Chapter~7]{Ouh} and references therein.\\

It follows from the definition of the operator $L$ and the fact that $U$ is a harmonic function of $H$ that $e^{-tL}{\bf 1} = {\bf 1}$. In other words,
$$\int_{{\bf R}^N} G(x,y,t)\omega(y) dy = 1.$$
This together with the  doubling property \eqref{eq:2.1} and the Gaussian upper bound \eqref{eq:4.9} 
imply the diagonal lower bound
\begin{equation}\label{eq:4.10}
G(x,x,t) \ge  \frac{C}{\omega(B(x,\sqrt{t}))}
\end{equation} 
for some constant $C > 0$. See,  e.g., \cite{Cou} and \cite[Chapter 7]{Ouh}.
Next, one extends this  diagonal lower bound to $x$ and $y$ near the diagonal. 
In order to do this one needs  
the H\"older continuity of the heat kernel $G(t,x,y)$. 
This latter property follows from the Harnack inequality. 
The  H\"older continuity  is also proved in 
\cite{CUR}, namely
\begin{equation}\label{eq:4.11}
| G(x,x,t) - G(x,y,t) | \le \frac{C t^{-\eta/2}} {\omega(B(x,\sqrt{t}))} | x- y|^\eta
\end{equation}
for some $\eta \in (0,1)$ and all $x, y$ and $t > 0$ such that $| x- y | \le \frac{1}{2} \sqrt{t}$. Using \eqref{eq:4.10} and \eqref{eq:4.11} one obtains easily
$$
G(x,y,t) \ge  \frac{C}{\sqrt{\omega(B(x,\sqrt{t}))}\sqrt{\omega(B(y,\sqrt{t}))}}
$$
for $x, y \in {\bf R}^N$ and $t > 0$ such that  $| x- y | \le \delta \sqrt{t}$ for some constant $\delta > 0$. Finally, the Gaussian lower bound
$$
G(x,y,t) \ge  \frac{C}{\sqrt{\omega(B(x,\sqrt{t}))}\sqrt{\omega(B(y,\sqrt{t}))}}\exp\left(-\frac{|x-y|^2}{Ct}\right)
$$
follows by a chain argument and the semigroup property. 
See again, e.g., \cite{Cou} and \cite[Chapter 7]{Ouh}.  
The equality \eqref{eq:4.1} gives the lower estimate of Theorem \ref{Theorem:1.2}. 
Thus the proof of Theorem \ref{Theorem:1.2} is complete.
\section{Non $A_2$~weight}\label{section:5}
In this section we study upper bounds of $p=p(x,y,t)$ 
without the assumption that $U^2$ is an $A_2$~weight on ${\bf R}^N$, 
and prove Theorem~\ref{Theorem:1.3}. 
In what follows, we set 
$$
d:=A^+(\lambda_2)
\quad\mbox{if $H$ is subcritical},
\qquad
d:=A^-(\lambda_2)
\quad\mbox{if $H$ is critical}.
$$
The first lemma follows from a similar argument as in the proof of Lemma~\ref{Lemma:4.1}. 
\begin{lemma}
\label{Lemma:5.1}
Let $V$ be a continuous function on $(0,\infty)$ satisfying \eqref{eq:1.2}.  
Assume that $H:=-\Delta+V(|x|)$ is nonnegative and let $U$ be as in Theorem~{\rm\ref{Theorem:1.1}}.
Then, for any $\epsilon\in(0,1]$, there exists a constant $C$ such that 
\begin{equation}
\label{eq:5.1}
p(x,y,t)\le Ct^{-\frac{N}{2}}\exp\left(-\frac{|x-y|^2}{Ct}\right)
\end{equation}
for all $x$, $y\in{\bf R}^N\setminus B(0,\epsilon\sqrt{t})$ and $t>0$. 
In particular, 
\begin{equation}
\label{eq:5.2}
p(x,y,t)\le Ct^{-\frac{N}{2}}\frac{U(\min\{|x|,\sqrt{t}\})U(\min\{|y|,\sqrt{t}\})}{U(\sqrt{t})^2}
\exp\left(-\frac{|x-y|^2}{Ct}\right)
\end{equation}
for all $x$, $y\in{\bf R}^N\setminus B(0,\epsilon\sqrt{t})$ and $t>0$. 
\end{lemma}
{\bf Proof.} 
Let $\epsilon\in(0,1]$, $x$, $y\in{\bf R}^N\setminus B(0,\epsilon\sqrt{t})$ and $t>0$. 
Similarly to the proof of Lemma~\ref{Lemma:4.1}, 
let $\lambda\in{\bf R}$ and 
let $\psi$ be a bounded smooth function $\psi$ on ${\bf R}^N$ such that $|\nabla\psi|\le 1$ on ${\bf R}^N$. 
For any $f_0\in L^2(B(y,\epsilon\sqrt{t}))$, set 
$$
f(\xi,s):=\int_{B(y,\epsilon\sqrt{t})}p(\xi,z,s)e^{-\lambda\psi(z)}f_0(z)\,dz,
\quad
F(\xi,s):=e^{\lambda\psi(\xi)}f(\xi,s). 
$$
Then it follows from the nonnegativity of $H$ that 
\begin{equation*}
\begin{split}
 & \frac{d}{ds}\int_{{\bf R}^N}F(\xi,s)^2\,d\xi\\
 & =-2\int_{{\bf R}^N}e^{2\lambda\psi}\left[2\lambda f\,\nabla\psi\cdot\nabla f
 +|\nabla f|^2+Vf^2\right]\,d\xi\\
 & =-2\int_{{\bf R}^N}\left[|\nabla(e^{\lambda\psi}f)|^2+V(e^{\lambda\psi}f)^2\right]\,d\xi
 +2\lambda^2\int_{{\bf R}^N}e^{2\lambda\psi}f^2|\nabla\psi|^2\,d\xi\\
 & \le 2\lambda^2\int_{{\bf R}^N}e^{2\lambda\psi}f^2|\nabla\psi|^2\,d\xi
 \le 2\lambda^2\int_{{\bf R}^N}F(\xi,s)^2\,d\xi,
\end{split}
\end{equation*}
which implies that 
$$
\int_{{\bf R}^N}F(\xi,s)^2\,d\xi\le e^{2\lambda^2 s}\int_{B(y,\sqrt{t})}f_0(\xi)^2\,d\xi. 
$$

Let $0<\eta\le t$. Set 
$$
\tilde{f}(\xi,s):=f(x+\delta\epsilon\xi,\eta+\delta^2\epsilon^2 s)\quad\mbox{with}\quad\delta=\sqrt{\eta}/4,
\qquad
\tilde{V}(\xi):=\delta^2\epsilon^2V(x+\delta\epsilon\xi). 
$$
Since $f$ satisfies $\partial_s f=\Delta_\xi f-V(\xi)f$ on ${\bf R}^N\times(0,\infty)$, 
we have
$$
\partial_s\tilde{f}=\Delta_\xi \tilde{f}-\tilde{V}(\xi)\tilde{f}
\quad\mbox{in}\quad{\bf R}^N\times\left(-16\epsilon^{-2},\infty\right).
$$
Furthermore, 
\begin{equation}
\label{eq:5.3}
|x+\delta\epsilon\xi|
\ge|x|-\frac{\epsilon}{4}\sqrt{\eta}|\xi|
\ge\epsilon\sqrt{t}-\frac{\epsilon}{4}\sqrt{t}|\xi|
\ge\frac{\epsilon}{2}\sqrt{t}\ge\frac{\epsilon}{2}\sqrt{\eta}
\end{equation}
for $x\in{\bf R}^N\setminus B(0,\epsilon\sqrt{t})$ and $\xi\in B(0,2)$. 
Since $|V(|x|)|\le C|x|^{-2}$ by \eqref{eq:1.2}, 
we deduce from \eqref{eq:5.3} that 
$$
|\tilde{V}(\xi)|\le C\delta^2\epsilon^2|x+\delta\epsilon\xi|^{-2}\le C,
\qquad
\xi\in B(0,2). 
$$
Then it follows from Proposition~\ref{Proposition:2.1} that 
$$
|\tilde{f}(0,0)|^2\le C\int_{-1}^0\int_{B(0,1)}|\tilde{f}|^2\,d\xi ds.
$$
Since $\eta-\delta^2\epsilon^2\ge 3\eta/4$ and $\delta\epsilon\le\epsilon\sqrt{t}$, 
it follows that 
\begin{equation*}
\begin{split}
f(x,\eta)^2
 & \le C(\delta\epsilon)^{-N-2}\int_{\eta-\delta^2\epsilon^2}^\eta\int_{B(x,\delta\epsilon)}|f(\xi,s)|^2\,d\xi ds\\
 & \le C\eta^{-\frac{N}{2}-1}\int_{3\eta/4}^\eta\int_{B(x,\epsilon\sqrt{t})}|f(\xi,s)|^2\,d\xi ds
\end{split}
\end{equation*}
for $x\in{\bf R}^N\setminus B(0,\epsilon\sqrt{t})$ and $0<\eta\le t\le 1$. 
Then we apply a similar argument as in the proof of Lemma~\ref{Lemma:4.1} 
to obtain \eqref{eq:5.1}. 
Furthermore, it follows from Theorem~\ref{Theorem:1.1} that 
\begin{equation}
\label{eq:5.4}
U(\min\{|x|,\epsilon\sqrt{t}\})\asymp U(\min\{|x|,\sqrt{t}\})
\quad\mbox{and}\quad 
U(\epsilon\sqrt{t})\asymp U(\sqrt{t})
\end{equation}
for $x\in{\bf R}^N$ and $t>0$. 
Then we deduce from \eqref{eq:5.1} and \eqref{eq:5.4} that 
\begin{equation*}
\begin{split}
 p(x,y,t) & \le Ct^{-\frac{N}{2}}\frac{U(\min\{|x|,\epsilon\sqrt{t}\})U(\min\{|y|,\epsilon\sqrt{t}\})}{U(\epsilon\sqrt{t})^2}
 \exp\left(-\frac{|x-y|^2}{Ct}\right)\\
  & \le Ct^{-\frac{N}{2}}\frac{U(\min\{|x|,\sqrt{t}\})U(\min\{|y|,\sqrt{t}\})}{U(\sqrt{t})^2}
 \exp\left(-\frac{|x-y|^2}{Ct}\right)
\end{split}
\end{equation*}
for all $x$, $y\in{\bf R}^N\setminus B(0,\epsilon\sqrt{t})$ and $t>0$. 
So we have \eqref{eq:5.2}, and the proof is complete. 
$\Box$\vspace{5pt}
\newline
Combining Lemma~\ref{Lemma:5.1} with Lemma~\ref{Lemma:4.1},  
we obtain upper estimates of $p=p(x,y,t)$ 
in the case where $0<t\le 1$ and $A^+(\lambda_1)<N/2$. 
\begin{lemma}
\label{Lemma:5.2}
Assume the same conditions as in Theorem~{\rm\ref{Theorem:1.3}} and $A^+(\lambda_1)<N/2$. 
Then there exists a constant $C$ such that 
\begin{equation}
\label{eq:5.5}
p(x,y,t)\le Ct^{-\frac{N}{2}}\frac{U(\min\{|x|,\sqrt{t}\})U(\min\{|y|,\sqrt{t}\})}{U(\sqrt{t})^2}
\exp\left(-\frac{|x-y|^2}{Ct}\right)
\end{equation}
for all $x$, $y\in{\bf R}^N$ and $0<t\le 1$. 
\end{lemma}
{\bf Proof.}
Let $G=G(x,y,t)$ be as in Section~4. 
Let $0<t\le 1$. The proof is divided into the following four cases: 
\begin{equation*}
\begin{array}{ll}
\mbox{(i)}\quad x,y\in{\bf R}^N\setminus B(0,\sqrt{t});
\qquad 
 & \mbox{(ii)}\quad
 x,y\in B(0,\sqrt{t});\vspace{5pt}\\
\mbox{(iii)}\quad
x\in{\bf R}^N\setminus B(0,\sqrt{t}),\,\,
y\in B(0,\sqrt{t});\quad
 & \mbox{(iv)}\quad
x\in B(0,\sqrt{t}),\,\,
y\in {\bf R}^N\setminus B(0,\sqrt{t}).
\end{array}
\end{equation*}
In case~(i) \eqref{eq:5.5} follows from Lemma~\ref{Lemma:5.1}. 
So we have only to consider cases~(ii), (iii) and (iv). 

Consider case (ii). 
It follows from Theorem~\ref{Theorem:1.1} that 
$U(|x|)\thicksim |x|^{A^+(\lambda_1)}$ as $|x|\to 0$. 
Combining \eqref{eq:1.4} with the assumption $A^+(\lambda_1)<N/2$, 
we see that $A^+(\lambda_1)\in(-N/2,N/2)$, which means that 
$\omega(x)=U(|x|)^2$ is an $A_2$~weight on  $B(0,2)$. 
Then Lemma~\ref{Lemma:4.1} implies that 
\begin{equation}
\label{eq:5.6}
G(x,y,t)\le 
\frac{C}{\sqrt{\omega(B(x,\sqrt{t}))}\sqrt{\omega(B(y,\sqrt{t}))}}\exp\left(-\frac{|x-y|^2}{Ct}\right).
\end{equation}
Furthermore, it follows from Theorem~\ref{Theorem:1.1} and \eqref{eq:2.1} that 
\begin{equation}
\label{eq:5.7}
\omega(B(\xi,\sqrt{s}))\asymp \omega(B(\xi,2\sqrt{s}))
\ge \omega(B(0,\sqrt{s}))\asymp s^{\frac{N}{2}+A^+(\lambda_1)}\asymp s^{\frac{N}{2}}U(\sqrt{s})^2
\end{equation}
for $\xi\in B(0,\sqrt{s})$ and $0<s\le 1$. 
By \eqref{eq:4.1}, \eqref{eq:5.6} and \eqref{eq:5.7} we obtain 
\begin{equation*}
\begin{split}
p(x,y,t)
 & \le\frac{CU(|x|)U(|y|)}{\sqrt{\omega(B(x,\sqrt{t}))}\sqrt{\omega(B(y,\sqrt{t}))}}\exp\left(-\frac{|x-y|^2}{Ct}\right)\\
  & \le Ct^{-\frac{N}{2}}\frac{U(|x|)U(|y|)}{U(\sqrt{t})^2}\exp\left(-\frac{|x-y|^2}{Ct}\right)\\
  & =Ct^{-\frac{N}{2}}\frac{U(\min\{|x|,\sqrt{t}\})U(\min\{|y|,\sqrt{t}\})}{U(\sqrt{t})^2}\exp\left(-\frac{|x-y|^2}{Ct}\right),
\end{split}
\end{equation*}
which implies \eqref{eq:5.5} in case~(ii).

Consider case~(iii). 
Set $\tilde{y}:=\sqrt{2t}y/|y|$ and $g(\xi,s):=G(x,\xi,s)$. 
Recalling that $w$ is $A_2$~weight on $B(0,2)$, 
we apply Lemma~\ref{Lemma:2.1} to $g$ to obtain  
$$
g(y,t)\le Cg(\tilde{y},2t)\exp\left(C\frac{|y-\tilde{y}|^2}{t}\right)
\le Cg(\tilde{y},2t),
$$
which together with \eqref{eq:4.1} implies 
\begin{equation}
\label{eq:5.8}
p(x,y,t)\le C\frac{U(|y|)}{U(|\tilde{y}|)}p(x,\tilde{y},2t)
=C\frac{U(|y|)}{U(\sqrt{2t})}p(x,\tilde{y},2t)
\le C\frac{U(|y|)}{U(\sqrt{t})}p(x,\tilde{y},2t).
\end{equation}
Since $|x|\ge\sqrt{t}=\epsilon\sqrt{2t}$ with $\epsilon=1/\sqrt{2}$, 
applying Lemma~\ref{Lemma:5.1}, we have 
\begin{equation*}
\begin{split}
p(x,\tilde{y},2t)
 & \le C(2t)^{-\frac{N}{2}}\frac{U(\min\{|x|,\sqrt{2t}\})U(\min\{|\tilde{y}|,\sqrt{2t}\})}{U(\sqrt{2t})^2}
\exp\left(-\frac{|x-\tilde{y}|^2}{2Ct}\right)\\
 & \le Ct^{-\frac{N}{2}}\frac{U(\min\{|x|,\sqrt{2t}\})}{U(\sqrt{2t})}
\exp\left(-\frac{|x|^2}{Ct}\right)\\
 & \le Ct^{-\frac{N}{2}}\frac{U(\min\{|x|,\sqrt{t}\})}{U(\sqrt{t})}
\exp\left(-\frac{|x-y|^2}{Ct}\right).
\end{split}
\end{equation*}
This together with \eqref{eq:5.8} implies \eqref{eq:5.5} in case~(iii). 
Since $p(x,y,t)=p(y,x,t)$, 
we also obtain \eqref{eq:5.5} in case~(iv). 
Thus Lemma~\ref{Lemma:5.2} follows.
$\Box$\vspace{5pt}

Next we obtain upper estimates of $p=p(x,y,t)$ 
in the case where $0<t\le 1$ and $A^+(\lambda_1)\ge N/2$. 
\begin{lemma}
\label{Lemma:5.3}
Assume the same conditions as in Theorem~{\rm\ref{Theorem:1.3}} and $A^+(\lambda_1)\ge N/2$. 
Then there exists a constant $C$ such that 
\begin{equation}
\label{eq:5.9}
p(x,y,t)\le Ct^{-\frac{N}{2}}\frac{U(\min\{|x|,\sqrt{t}\})U(\min\{|y|,\sqrt{t}\})}{U(\sqrt{t})^2}
\exp\left(-\frac{|x-y|^2}{Ct}\right)
\end{equation}
for all $x$, $y\in{\bf R}^N$ and $0<t<1$. 
\end{lemma}
For this aim, we prepare the following lemma,  
which is useful to obtain upper estimates of $p=p(x,y,t)$ inside a parabolic cone. 
A similar lemma has been used in the study of the behavior of the solutions of 
the heat equation with a potential (see e.g., \cite{IIY01, IIY02, IK01, IK04}). 
\begin{lemma}
\label{Lemma:5.4}
Assume the same conditions as in Theorem~{\rm\ref{Theorem:1.1}}. 
Let $T\ge 0$. 
Furthermore, assume that 
$$
\zeta(t):=
t^{\gamma_1}[\log(c+t)]^{\gamma_2}
$$
is monotone decreasing on $(T,\infty)$, 
where $\gamma_1$, $\gamma_2\in{\bf R}$ and $c>1$. 
Let $\kappa>0$ be such that 
\begin{equation}
\label{eq:5.10}
-s\zeta'(s)\le\kappa\zeta(s),\qquad s\in(T,\infty). 
\end{equation}
Define 
\begin{equation*}
\begin{split}
F[U](x):= & U(|x|)\int_0^{|x|} s^{1-N}[U(s)]^{-2}
\left(\int_0^s \tau^{N-1}U(\tau)^2\,d\tau\right)\,ds,\\
w(x,t):= & \zeta(s)\left[U(|x|)-\kappa s^{-1}F[U](x)\right].
\end{split}
\end{equation*}
Then
$$
\partial_t w\ge\Delta w-V(|x|)w
\quad\mbox{in}\quad{\bf R}^N\times(T,\infty).
$$
\end{lemma}
{\bf Proof.}
It follows that 
$\Delta F-V(|x|)F=U(|x|)$ for $x\in{\bf R}^N$. 
This together with \eqref{eq:5.10} implies 
$$
\partial_t w-\Delta w+V(|x|)w
\ge [\zeta'(t)+\kappa t^{-1}\zeta(t)]U(x)\ge 0,
\quad x\in{\bf R}^N,\,t\in(T,\infty). 
$$
Thus Lemma~\ref{Lemma:5.4} follows. 
$\Box$\vspace{5pt}
\newline
{\bf Proof of Lemma~\ref{Lemma:5.3}.}
For any $\sigma>0$, we define 
$$
H_\sigma:=-\Delta+V_\sigma(|x|),
\qquad
V_\sigma(|x|):=\frac{V(|x|)}{U(|x|)+\sigma}U(|x|).
$$
Let $p_\sigma=p_\sigma(x,y,t)$ be the fundamental solution corresponding to $e^{-t H_\sigma}$. 
It follows from Theorem~\ref{Theorem:1.1} and \eqref{eq:1.2} that 
\begin{equation}
\label{eq:5.11}
|V_\sigma(|x|)|\le |V(|x|)|\le C|x|^{-2}\,\,\,\mbox{in}\,\,\,{\bf R}^N,
\quad
V_\sigma(r)\thicksim \lambda_1r^{-2+A^+(\lambda_1)}
\quad\mbox{as}\quad r\to 0. 
\end{equation}
In particular, since $A^+(\lambda_1)>0$, we see that
$V_\sigma\in L^q_{{\rm loc}}({\bf R}^N)$ for some $q>N/2$. 
Furthermore, $U_\sigma:=U+\sigma$ is a positive harmonic function for $H_\sigma$
and 
$$
\int_{{\bf R}^N}\left[|\nabla\varphi|^2+V_\sigma\varphi^2\right]\,dz
=\int_{{\bf R}^N}\left|\nabla\frac{\varphi}{U_\sigma}\right|^2U_\sigma^2\,dz
\ge 0
$$
for all $\varphi\in C_0^\infty({\bf R}^N)$, which means that $H_\sigma$ is nonnegative on $L^2({\bf R}^N)$. 

In the proof, the letter $C_*$ denotes a generic constant independent of $x$, $y$, $t$ and $\sigma$. 
Since $H_\sigma$ is nonnegative, for any $\epsilon\in(0,1]$, 
we apply Lemma~\ref{Lemma:5.1} with the aid of \eqref{eq:5.11} to obtain 
\begin{equation}
\label{eq:5.12}
p_\sigma(t,x,y)\le C_*t^{-\frac{N}{2}}\exp\left(-\frac{|x-y|^2}{C_*t}\right)
\end{equation}
for all $x$, $y\in{\bf R}^N\setminus B(0,\epsilon\sqrt{t})$ and $t>0$. 
On the other hand, 
since $U_\sigma^2$ is an $A_2$~weight on $B(0,R)$ for any $R>0$, 
we apply a similar argument as in the proof of Lemma~\ref{Lemma:5.2} 
to obtain 
\begin{equation}
\label{eq:5.13}
p_\sigma(x,y,t)
\le C_{R,\sigma}
t^{-\frac{N}{2}}\frac{U_\sigma(\min\{|x|,\sqrt{t}\})U_\sigma(\min\{|y|,\sqrt{t}\})}{U_\sigma(\sqrt{t})^2}
\exp\left(-\frac{|x-y|^2}{C_{R,\sigma}t}\right)
\end{equation}
for all $x$, $y\in B(0,R)$ and $0<t\le 1$, 
where $C_{R,\sigma}$ is a constant depending on $R$ and $\sigma$. 

Let $\epsilon$ be a sufficiently small positive constant to be chosen later. 
Let $x$, $y\in{\bf R}^N$ and $0<t\le 1$. 
In what follows, we divide the proof into the following four cases: 
\begin{equation*}
\begin{array}{ll}
\mbox{(i)}\quad x,y\in{\bf R}^N\setminus B(0,\epsilon\sqrt{t});
\qquad 
 & \mbox{(ii)}\quad
 x,y\in B(0,\epsilon\sqrt{t});\vspace{5pt}\\
\mbox{(iii)}\quad
x\in{\bf R}^N\setminus B(0,\epsilon\sqrt{t}),\,\,
y\in B(0,\epsilon\sqrt{t});\quad
 & \mbox{(iv)}\quad
x\in B(0,\epsilon\sqrt{t}),\,\,
y\in {\bf R}^N\setminus B(0,\epsilon\sqrt{t}).
\end{array}
\end{equation*}
Similarly to Lemma~\ref{Lemma:5.2}, 
by Lemma~\ref{Lemma:5.1} we have \eqref{eq:5.9} in case~(i). 

We consider case~(iii). 
Define 
\begin{equation*}
\begin{split}
D_\epsilon(t):= & \left\{(\xi,s)\in{\bf R}^N\times(0,t]\,:\,|\xi|<\epsilon\sqrt{s}\right\},\\
\partial_pD_\epsilon(t):= & \left\{(\xi,s)\in{\bf R}^N\times[0,t]\,:\,|\xi|=\epsilon\sqrt{s}\right\}.
\end{split}
\end{equation*}
Let $\kappa:=(N+A^+(\lambda_1))/2$ and set 
\begin{equation}
\label{eq:5.14}
\begin{split}
 & v(\xi,s):=p_\sigma(x,\xi,s),
\quad
w(\xi,s):=s^{-\frac{N+A^+(\lambda_1)}{2}}\left[U_\sigma(|\xi|)-\kappa s^{-1}F[U_\sigma](|\xi|)\right],\\
 & z(\xi,s):=v(\xi,s)-\gamma\exp\left(-\frac{|x|^2}{\gamma t}\right)w(\xi,s),
 \end{split}
\end{equation}
where $\gamma$ is a positive constant. 
It follows from Lemma~\ref{Lemma:5.4} that  
\begin{equation}
\label{eq:5.15}
\partial_s z\le\Delta z-V_\sigma(|\xi|)z
\quad\mbox{in}\quad {\bf R}^N\times(0,\infty). 
\end{equation}
Since 
$U(r)\asymp r^{A^+(\lambda_1)}$ on $(0,1)$ and $A^+(\lambda_1)>0$, 
we have 
\begin{equation}
\label{eq:5.16}
\begin{split}
 & F[U_\sigma](|\xi|)=U_\sigma(|\xi|)\int_0^{|\xi|} s^{1-N}[U_\sigma(s)]^{-2}
\left(\int_0^s \tau^{N-1}U_\sigma(\tau)^2\,d\tau\right)\,ds\\
 & \qquad
\le U_\sigma(|\xi|)\int_0^{|\xi|} s^{1-N}[C_*^{-1}s^{A^+(\lambda_1)}+\sigma]^{-2}
\left(\int_0^s \tau^{N-1}(C_*r^{A^+(\lambda_1)}+\sigma)^2\,d\tau\right)\,ds\\
 & \qquad
\le U_\sigma(|\xi|)\int_0^{|\xi|} s^{1-N}[C_*^{-1}\sigma^{-1}s^{A^+(\lambda_1)}+1]^{-2}\cdot
\frac{1}{N}s^N(C_*\sigma^{-1}s^{A^+(\lambda_1)}+1)^2\,ds\\
 & \qquad
\le\frac{C_*}{2N}|\xi|^2U_\sigma(|\xi|)
\le\frac {C_*\epsilon}{2N}sU_\sigma(|\xi|),
\qquad (\xi,s)\in D_\epsilon(t). 
\end{split}
\end{equation}
Taking a sufficiently small $\epsilon>0$ if necessary,  
by \eqref{eq:5.14} and \eqref{eq:5.16} we obtain 
$$
w(\xi,s)\ge\frac{1}{2}s^{-\frac{N+A^+(\lambda_1)}{2}}U_\sigma(|\xi|)
=\frac{1}{2}s^{-\frac{N+A^+(\lambda_1)}{2}}[U(|\xi|)+\sigma]
\quad\mbox{in}\quad D_\epsilon(t).
$$
This implies that 
\begin{equation}
\label{eq:5.17}
w(\xi,s)\ge\frac{1}{2}s^{-\frac{N+A^+(\lambda_1)}{2}}\sigma
\quad\mbox{in}\quad D_\epsilon(t),
\qquad\qquad\quad
\end{equation}
\begin{equation}
\label{eq:5.18}
\begin{split}
\qquad\qquad\qquad\,\,
w(\xi,s) & \ge\frac{1}{2}s^{-\frac{N+A^+(\lambda_1)}{2}}U(|\xi|)
\ge C_*s^{-\frac{N+A^+(\lambda_1)}{2}}(\epsilon\sqrt{s})^{A^+(\lambda_1)}\\
 & \ge C_*\epsilon^{A^+(\lambda_1)}s^{-\frac{N}{2}}
\quad\mbox{on}\quad \partial_p D_\epsilon(t)\setminus\{(0,0)\}.
\end{split}
\end{equation}
On the other hand, it follows from \eqref{eq:5.12} that
\begin{equation}
\label{eq:5.19}
v(\xi,s)=p_\sigma(x,\xi,s)
\le C_*s^{-\frac{N}{2}}\exp\left(-\frac{|x-\xi|^2}{C_*s}\right)
\le C_*s^{-\frac{N}{2}}\exp\left(-\frac{|x|^2}{C_*t}\right)
\end{equation}
on $ \partial_p D_\epsilon(t)\setminus\{(0,0)\}$.
Then, by \eqref{eq:5.18} and \eqref{eq:5.19}, taking a sufficiently large constant $\gamma$ if necessary, 
we have  
\begin{equation}
\label{eq:5.20}
z(\xi,s)\le 0\quad\mbox{on}\quad \partial_p D_\epsilon(t)\setminus\{(0,0)\}.
\end{equation}
On the other hand, since $|x|\ge\epsilon\sqrt{t}$ and $A^+(\lambda_1)>0$, 
by \eqref{eq:5.13} we see that 
$$
\lim_{s\to 0} s^{\frac{N+A^+(\lambda_1)}{2}}v(\xi,s)=0
$$
uniformly for $\xi$ in a neighborhood of the origin. 
Then, by \eqref{eq:5.17} we see that 
\begin{equation}
\label{eq:5.21}
z(\xi,s)\le 0
\end{equation}
for $(\xi,s)\in D_\epsilon(t)$ if $s$ is sufficiently small. 
Therefore, by \eqref{eq:5.15}, \eqref{eq:5.20} and \eqref{eq:5.21} 
we apply the comparison principle to obtain 
$z\le 0$ on $D_\epsilon(t)$. 
This together with \eqref{eq:5.14} implies that 
$$
p_\sigma(x,\xi,s)=v(\xi,s)\le \gamma\exp\left(-\frac{|x|^2}{\gamma t}\right)w(\xi,s)
\le\gamma s^{-\frac{N+A^+(\lambda_1)}{2}}U_\sigma(|\xi|)\exp\left(-\frac{|x|^2}{\gamma t}\right)
$$
for $(\xi,s)\in D_\epsilon(t)$ and $0<\sigma\le 1$.  
Taking $(\xi,s)=(y,t)$, we obtain  
\begin{equation}
\label{eq:5.22}
\begin{split}
p_\sigma(x,y,t) & \le \gamma t^{-\frac{N+A^+(\lambda_1)}{2}}U_\sigma(|y|)\exp\left(-\frac{|x|^2}{\gamma t}\right)\\
 & \le C_*t^{-\frac{N+A^+(\lambda_1)}{2}}U_\sigma(|y|)\exp\left(-\frac{|x-y|^2}{C_*t}\right)\\
 & \le C_*t^{-\frac{N}{2}}\frac{U(\min\{|y|,\epsilon\sqrt{t}\})+\sigma}{U(\sqrt{t})}\exp\left(-\frac{|x-y|^2}{C_*t}\right). 
\end{split}
\end{equation}
Passing to the limit as $\sigma\to 0$, 
we deduce that 
\begin{equation*}
\begin{split}
p(x,y,t)
 & \le C_*t^{-\frac{N}{2}}\frac{U(\min\{|y|,\epsilon\sqrt{t}\})}{U(\sqrt{t})}\exp\left(-\frac{|x-y|^2}{C_*t}\right)\\
 & \le C_*t^{-\frac{N}{2}}\frac{U(\min\{|y|,\sqrt{t}\})}{U(\sqrt{t})}\exp\left(-\frac{|x-y|^2}{C_*t}\right),
\end{split}
\end{equation*}
which means that \eqref{eq:5.9} holds in case~(iii). 
Since $p(x,y,t)=p(y,x,t)$, 
we also obtain \eqref{eq:5.9} in case~(iv).

It remains to consider case~(ii).
Let $\tilde{\kappa}:=(N+2A^+(\lambda_1))/2$ and set
\begin{equation}
\label{eq:5.23}
\begin{split}
 & \tilde{v}(\xi,s):=p_\sigma(\xi,y,s),
\quad
\tilde{w}(\xi,s):=s^{-\frac{N+2A^+(\lambda_1)}{2}}\left[U_\sigma(|\xi|)-\tilde{\kappa} s^{-1}F[U_\sigma](|\xi|)\right],\\
 & \tilde{z}(\xi,s):=\tilde{v}(\xi,s)-\gamma'
[\gamma' U(\min\{|y|,\epsilon\sqrt{t}\})+\sigma]\exp\left(-\frac{|y|^2}{\gamma' t}\right)\tilde{w}(\xi,s),
\end{split}
\end{equation}
where $\gamma'$ is a positive constant. 
It follows from Lemma~\ref{Lemma:5.4} that 
\begin{equation}
\label{eq:5.24}
\partial_s \tilde{z}\le\Delta_\xi \tilde{z}-V_\sigma(|\xi|)\tilde{z}
\quad\mbox{in}\quad {\bf R}^N\times(0,\infty). 
\end{equation}
For $(\xi,s)\in\partial_p D_\epsilon(t)\setminus\{(0,0)\}$, 
we see that 
$\xi\in{\bf R}^N\setminus B(0,\epsilon\sqrt{s})$. 
Since 
$U(r)\asymp r^{A^+(\lambda_1)}$ on $(0,1)$ and $A^+(\lambda_1)>0$, 
we apply \eqref{eq:5.9} and \eqref{eq:5.22} to obtain 
\begin{equation}
\label{eq:5.25}
\begin{split}
\tilde{v}(\xi,s) & \le C_*s^{-\frac{N}{2}}
\frac{U(\min\{|y|,\epsilon\sqrt{s}\})+\sigma}{U(\sqrt{s})}\exp\left(-\frac{|\xi-y|^2}{C_*s}\right)\\
 & \le C_*s^{-\frac{N+A^+(\lambda_1)}{2}}[C_*U(\min\{|y|,\epsilon\sqrt{t}\})+\sigma]\exp\left(-\frac{|y|^2}{C_*s}\right)\\
 & \le C_*s^{-\frac{N+A^+(\lambda_1)}{2}}[C_*U(\min\{|y|,\epsilon\sqrt{t}\})+\sigma]\exp\left(-\frac{|y|^2}{C_*t}\right)
\end{split}
\end{equation}
for $(\xi,s)\in\partial_p D_\epsilon(t)\setminus\{(0,0)\}$. 
On the other hand, 
taking a sufficiently small $\epsilon>0$ if necessary, 
by \eqref{eq:5.16} we have 
$$
\tilde{w}(\xi,s)\ge\frac{1}{2}s^{-\frac{N+2A^+(\lambda_1)}{2}}U_\sigma(|x|)\quad\mbox{in}\quad D_\epsilon(t).
$$
Then, similarly to \eqref{eq:5.17} and \eqref{eq:5.18}, 
we see that   
\begin{equation}
\label{eq:5.26}
\tilde{w}(\xi,s)\ge\frac{1}{2}s^{-\frac{N+2A^+(\lambda_1)}{2}}\sigma
\quad\mbox{in}\quad D_\epsilon(t),
\qquad\qquad\quad
\end{equation}
\begin{equation}
\label{eq:5.27}
\begin{split}
\tilde{w}(\xi,s) & \ge\frac{1}{2}s^{-\frac{N+2A^+(\lambda_1)}{2}}U(|\xi|)
\ge C_*s^{-\frac{N+2A^+(\lambda_1)}{2}}(\epsilon\sqrt{s})^{A^+(\lambda_1)}\\
 & \ge C_*\epsilon^{A^+(\lambda_1)}s^{-\frac{N+A^+(\lambda_1)}{2}}
\quad\mbox{on}\quad \partial D_\epsilon(t).
\end{split}
\end{equation}
Taking a sufficiently large constant $\gamma'$ if necessary, 
by \eqref{eq:5.25} and \eqref{eq:5.27} 
we have  
\begin{equation}
\label{eq:5.28}
\tilde{z}(\xi,s)\le 0\quad\mbox{on}\quad\partial_p D_\epsilon(t)\setminus\{(0,0)\}. 
\end{equation}
Furthermore,  
by \eqref{eq:5.13} we see that 
$$
\lim_{s\to 0}s^{\frac{N+2A^+(\lambda_1)}{2}}\tilde{v}(\xi,s)=0
$$
uniformly for $\xi$ in a neighborhood of the origin. 
This together with \eqref{eq:5.17} implies that 
\begin{equation}
\label{eq:5.29}
\tilde{z}(\xi,s)\le 0 
\end{equation}
for $(\xi,s)\in D_\epsilon(t)$ if $s$ is sufficiently small. 
Therefore, 
by \eqref{eq:5.24}, \eqref{eq:5.28} and \eqref{eq:5.29} 
we apply the comparison principle to obtain 
$\tilde{z}\le 0$ on $D_\epsilon(t)$. 
This together with \eqref{eq:5.23} implies that 
\begin{equation*}
\begin{split}
p_\sigma(\xi,y,s)=\tilde{v}(\xi,s)
 & \le\gamma'
[\gamma' U(\min\{|y|,\epsilon\sqrt{t}\})+\sigma]\exp\left(-\frac{|y|^2}{\gamma' t}\right)\tilde{w}(\xi,s)\\
 & \le C_*s^{-\frac{N+2A^+(\lambda_1)}{2}}[C_*U(\min\{|y|,\epsilon\sqrt{t}\})+\sigma]U_\sigma(|y|)\exp\left(-\frac{|y|^2}{C_*t}\right)
\end{split}
\end{equation*}
for $(\xi,s)\in D_\epsilon(t)$. 
Taking $(\xi,s)=(x,t)$ and passing to the limit as $\sigma\to 0$, 
by \eqref{eq:5.4} we obtain  
\begin{equation*}
\begin{split}
p(x,y,t)=\lim_{\sigma\to 0}p_\sigma(x,y,t) & \le C_*{t}^{-\frac{N+2A^+(\lambda_1)}{2}}U(\min\{|y|,\epsilon\sqrt{t}\})
U(|x|)\exp\left(-\frac{|y|^2}{C_*t}\right)\\
 & \le C_*t^{-\frac{N+2A^+(\lambda_1)}{2}}U(\min\{|y|,\sqrt{t}\})U(|x|)\exp\left(-\frac{|x-y|^2}{C_*t}\right)\\
 & \le C_*t^{-\frac{N}{2}}\frac{U(\min\{|y|,\sqrt{t}\})U(\min\{|x|,\sqrt{t}\})}{U(\sqrt{t})^2}
 \exp\left(-\frac{|x-y|^2}{C_*t}\right). 
\end{split}
\end{equation*}
which means that \eqref{eq:5.9} holds in case~(ii). 
Thus Lemma~\ref{Lemma:5.2} follows. 
$\Box$\vspace{5pt}

We complete the proof of Theorem~\ref{Theorem:1.3}.
\vspace{5pt}
\newline
{\bf Proof of Theorem~\ref{Theorem:1.3}.}
Let $\epsilon$ be a sufficiently small positive constant. 
Due to Lemmas~\ref{Lemma:5.2} and \ref{Lemma:5.3}, 
it suffices to prove \eqref{eq:1.6} in the case $t>1$. 

Let $t>1$. 
Similarly to Lemma~\ref{Lemma:5.3}, 
the proof is divided into the following four cases: 
\begin{equation*}
\begin{array}{ll}
\mbox{(i)}\quad x,y\in{\bf R}^N\setminus B(0,\epsilon\sqrt{t});
\qquad 
 & \mbox{(ii)}\quad
 x,y\in B(0,\epsilon\sqrt{t});\vspace{5pt}\\
\mbox{(iii)}\quad
x\in{\bf R}^N\setminus B(0,\epsilon\sqrt{t}),\,
y\in B(0,\epsilon\sqrt{t});\quad
 & \mbox{(iv)}\quad
x\in B(0,\epsilon\sqrt{t}),\,
y\in {\bf R}^N\setminus B(0,\epsilon\sqrt{t}).
\end{array}
\end{equation*}
In case~(i), 
by Lemma~\ref{Lemma:5.1} we have \eqref{eq:1.6}.  

Consider case~(iii). 
Define 
\begin{equation*}
\begin{split}
E_\epsilon(t):= & \left\{(\xi,s)\in{\bf R}^N\times(1,t]\,:\,|\xi|<\epsilon\sqrt{s}\right\},\\
\partial_pE_\epsilon(t):= & \left\{(\xi,s)\in{\bf R}^N\times(1,t]\,:\,|\xi|=\epsilon\sqrt{s}\right\}\,\cup\,
\{(\xi,1)\in{\bf R}^N\times\{1\}\,:\,|\xi|\le\epsilon\}. 
\end{split}
\end{equation*}
Let 
$$
\zeta(s):=
\left\{
\begin{array}{ll}
s^{-\frac{N+d}{2}}[\log(2+s)]^{-1} & \mbox{if $\lambda_2=\lambda_*$ and $H$ is subcritical},\vspace{3pt}\\
s^{-\frac{N+d}{2}} & \mbox{otherwise}.
\end{array}
\right.
$$
It follows from Theorem~\ref{Theorem:1.1} that 
\begin{equation}
\label{eq:5.30}
\zeta(s)\asymp s^{-\frac{N}{2}}U(\sqrt{s})^{-1}\quad\mbox{in}\quad(1,\infty).
\end{equation}
Since $N+d>0$, 
we can find $\kappa>0$ such that $-s\zeta'(s)\le\kappa\zeta(s)$ on $(1,\infty)$. 
Set 
\begin{equation}
\label{eq:5.31}
\begin{split}
 & v(\xi,s):=p(x,\xi,s),
\quad
w(\xi,s):=\zeta(s)\left[U(|\xi|)-\kappa s^{-1}F[U](|\xi|)\right],\\
 & z(\xi,s):=v(\xi,s)-C_1\exp\left(-\frac{|x|^2}{C_1t}\right)w(\xi,s),
\end{split}
\end{equation}
where $C_1$ is a positive constant to be chosen later. 
It follows from Lemma~\ref{Lemma:5.4} that 
\begin{equation}
\label{eq:5.32}
\partial_sz\le\Delta_\xi z-V(|\xi|)z\quad\mbox{in}\quad {\bf R}^N\times(1,\infty).
\end{equation}
Since \eqref{eq:1.6} holds in the case $0<t\le 1$, 
we see that 
\begin{equation}
\label{eq:5.33}
\begin{split}
v(\xi,1)=p(x,\xi,1) & \le C\frac{U(\min\{|x|,1\})U(\min\{|\xi|,1\})}{U(1)^2}\exp\left(-\frac{|x-\xi|^2}{C}\right)\\
 & \le CU(|\xi|)\exp\left(-\frac{|x|^2}{Ct}\right)
\end{split}
\end{equation}
for $\xi\in B(0,1)$. 
Furthermore, by Lemma~\ref{Lemma:5.1} we have 
\begin{equation}
\label{eq:5.34}
v(\xi,s)\le Cs^{-\frac{N}{2}}\exp\left(-\frac{|x-\xi|^2}{Cs}\right)
\le Cs^{-\frac{N}{2}}\exp\left(-\frac{|x|^2}{Ct}\right)
\end{equation}
for $(\xi,s)\in{\bf R}^N\times(1,t)$ with $|{\xi}|=\epsilon\sqrt{s}$. 
On the other hand, 
taking a sufficiently small $\epsilon>0$ if necessary, 
by \eqref{eq:5.16} with $\sigma=0$ we have
$$
w(\xi,s)\ge\frac{1}{2}\zeta(s)U(|\xi|)\quad\mbox{in}\quad E_\epsilon(t). 
$$
In particular, 
\begin{equation}
\label{eq:5.35}
w(\xi,1)\ge\frac{1}{2}\zeta(1)U(|\xi|)
\end{equation}
for $\xi\in B(0,\epsilon)$. 
In addition, by \eqref{eq:5.30} we see that 
\begin{equation}
\label{eq:5.36}
w(\xi,s)\ge C^{-1}s^{-\frac{N}{2}}
\end{equation}
for $(\xi,s)\in{\bf R}^N\times(1,t)$ with $|\xi|=\epsilon\sqrt{s}$. 
Taking a sufficiently large $C_1$ if necessary, 
by \eqref{eq:5.33}, \eqref{eq:5.34}, \eqref{eq:5.35} and \eqref{eq:5.36} we have 
\begin{equation}
\label{eq:5.37}
z(\xi,s)\le 0\quad\mbox{on}\quad \partial_pE_\epsilon(t).
\end{equation}
Therefore, by \eqref{eq:5.32} and \eqref{eq:5.37} 
we apply the comparison principle to obtain 
$z\le 0$ on $E_\epsilon(t)$. 
This implies that 
$$
p(x,\xi,s) =v(\xi,s)\le C\exp\left(-\frac{|x|^2}{Ct}\right)w(\xi,t)
\le C\zeta(s)U(|\xi|)\exp\left(-\frac{|x|^2}{Ct}\right)
$$
on $E_\epsilon(t)$. 
Taking $(\xi,s)=(y,t)$, by \eqref{eq:5.4} and \eqref{eq:5.30} we obtain 
\begin{equation*}
\begin{split}
p(x,y,t) & \le C\zeta(t)U(|y|)\exp\left(-\frac{|x|^2}{Ct}\right)
\le C\zeta(t)U(|y|)\exp\left(-\frac{|x-y|^2}{Ct}\right)\\
 & \le Ct^{-\frac{N}{2}}\frac{U(\min\{|y|,\epsilon\sqrt{t}\})}{U(\sqrt{t})}\frac{U(\min\{|x|,\epsilon\sqrt{t}\})}{U(\epsilon\sqrt{t})}
 \exp\left(-\frac{|x-y|^2}{Ct}\right)\\
 & \le Ct^{-\frac{N}{2}}\frac{U(\min\{|x|,\sqrt{t}\})U(\min\{|y|,\sqrt{t}\})}{U(\sqrt{t})^2}\exp\left(-\frac{|x-y|^2}{Ct}\right). 
\end{split}
\end{equation*}
Thus \eqref{eq:1.6} holds in case~(iii). 
Since $p(x,y,t)=p(y,x,t)$, 
\eqref{eq:1.6} also holds in case~(iv).

It remains to prove \eqref{eq:1.6} in case~(ii). 
Set 
$$
S:=1\quad\mbox{if}\quad |x|\le 1,
\qquad
S:=|x|^2\quad\mbox{if}\quad |x|>1. 
$$
Then it follows that 
\begin{equation}
\label{eq:5.38}
1\le S<t,
\qquad
U(\min\{|x|,\sqrt{S}\})=U(|x|).  
\end{equation}
We show that 
\begin{equation}
\label{eq:5.39}
v(\xi,S)=p(x,\xi,S)\le
CS^{-\frac{N}{2}}U(\sqrt{S})^{-2}U(|x|)U(|\xi|)\exp\left(-\frac{|x|^2}{Ct}\right)
\end{equation}
for all $\xi\in B(0,\epsilon\sqrt{S})$. 
In the case $S=1$, that is $|x|\le 1$, 
combining Lemmas~\ref{Lemma:5.2} and \ref{Lemma:5.3} with \eqref{eq:5.38}, 
we have \eqref{eq:5.39}. 
So we consider the case $S>1$, that is $|x|>1$. 
Let $w$ and $z$ be as in \eqref{eq:5.31}. 
Then $z$ satisfies \eqref{eq:5.32} on ${\bf R}^N\times(1,S]$. 
Furthermore, by \eqref{eq:1.6} in cases~(i) and (iii) 
we see that 
\begin{equation}
\label{eq:5.40}
\begin{split}
v(\xi,s)=p(x,\xi,s)
 & \le Cs^{-\frac{N}{2}}\frac{U(\min\{|x|,\sqrt{s}\})
 U(\min\{|\xi|,\sqrt{s}\})}{U(\sqrt{s})^2}\exp\left(-\frac{|x-\xi|^2}{Cs}\right)\\
 & \le Cs^{-\frac{N}{2}}\exp\left(-\frac{|x|^2}{Ct}\right)
\end{split}
\end{equation}
for $(\xi,s)\in{\bf R}^N\times(1,S]$ with $|\xi|=\epsilon\sqrt{s}$ 
and that 
\begin{equation}
\label{eq:5.41}
\begin{split}
v(\xi,1)
 & \le C\frac{U(\min\{|x|,1\})U(\min\{|\xi|,1\})}{U(1)^2}\exp\left(-\frac{|x-\xi|^2}{C}\right)\\
 & \le CU(1)U(|\xi|)\exp\left(-\frac{|x|^2}{Ct}\right)
\end{split}
\end{equation}
for $\xi\in B(0,\epsilon)$. 
Then, 
by \eqref{eq:5.35}, \eqref{eq:5.36}, \eqref{eq:5.40} and \eqref{eq:5.41}, 
taking a sufficiently large $C_1$ if necessary, 
we see that 
$$
z(\xi,s)\le 0\quad\mbox{on}\quad \partial_pE_\epsilon(S).
$$
Then, by the comparison principle we see that $z(\xi,s)\le 0$ on $E_\epsilon(S)$. 
This together with  \eqref{eq:5.30} implies that 
\begin{equation*}
\begin{split}
v(\xi,S) & \le C_1S^{-\frac{N}{2}}U(\sqrt{S})^{-1}U(|\xi|)\exp\left(-\frac{|x|^2}{C_1t}\right)\\
 & =C_1S^{-\frac{N}{2}}U(\sqrt{S})^{-2}U(|x|)U(|\xi|)\exp\left(-\frac{|x|^2}{C_1t}\right)
\end{split}
\end{equation*}
for all $\xi\in B(0,\epsilon\sqrt{S})$,  
which implies \eqref{eq:5.39} in the case $S>1$. 
Therefore inequality~\eqref{eq:5.39} holds. 

We complete the proof of \eqref{eq:1.6} in case~(ii). 
Let 
$$
\tilde{\zeta}(s):=
\left\{
\begin{array}{ll}
s^{-\frac{N+2d}{2}}[\log(2+s)]^{-2} & \mbox{if $\lambda_2=\lambda_*$ and $H$ is subcritical},\vspace{3pt}\\
s^{-\frac{N+2d}{2}} & \mbox{otherwise}.
\end{array}
\right.
$$
It follows from Theorem~\ref{Theorem:1.1} and \eqref{eq:5.30} that 
\begin{equation}
\label{eq:5.42}
\tilde{\zeta}(s)\asymp s^{-\frac{N}{2}}U(\sqrt{s})^{-2}\asymp U(\sqrt{s})^{-1}\zeta(s)\quad\mbox{in}\quad(1,\infty).
\end{equation}
Since $N+2d>0$, 
we can find $\tilde{\kappa}>0$ such that $-s\tilde{\zeta}'(s)\le\tilde{\kappa}\tilde{\zeta}(s)$ on $(1,\infty)$. 
Set 
\begin{equation*}
\begin{split}
 & \tilde{w}(\xi,s):=\tilde{\zeta}(s)\left[U(|\xi|)-\tilde{\kappa} s^{-1}F[U](|\xi|)\right],\\
 & \tilde{z}(\xi,s):=v(\xi,s)-C_2U(|x|)\exp\left(-\frac{|x|^2}{C_2t}\right)\tilde{w}(\xi,s), 
\end{split}
\end{equation*}
where $C_2$ is a positive constant to be chosen later. 
Then, by Lemma~\ref{Lemma:5.4} we see that 
\begin{equation}
\label{eq:5.43}
\partial_s\tilde{z}\le\Delta_\xi\tilde{z}-V(|\xi|)\tilde{z}\quad\mbox{in}\quad {\bf R}^N\times(1,\infty).
\end{equation}
Since \eqref{eq:1.6} holds in case~(iv), 
it follows from \eqref{eq:5.30} and \eqref{eq:5.38} that 
\begin{equation}
\label{eq:5.44}
\begin{split}
v(\xi,s)=p(x,\xi,s)
 & \le Cs^{-\frac{N}{2}}\frac{U(\min\{|x|,\sqrt{s}\})
 U(\min\{|\xi|,\sqrt{s}\})}{U(\sqrt{s})^2}\exp\left(-\frac{|x-\xi|^2}{Cs}\right)\\
 & \le C\zeta(s)U(|x|)\exp\left(-\frac{|x|^2}{Ct}\right)
\end{split}
\end{equation}
for $(\xi,s)\in{\bf R}^N\times[S,t]$ with $|\xi|=\epsilon\sqrt{s}$. 
Furthermore, similarly to \eqref{eq:5.35} and \eqref{eq:5.36}, 
taking a sufficiently small $\epsilon>0$ if necessary, 
we deduce from \eqref{eq:5.42} that
\begin{equation}
\label{eq:5.45}
\tilde{w}(\xi,s)\ge \frac{1}{2}\tilde{\zeta}(s)U(|\xi|)\ge C^{-1}\zeta(s)
\end{equation}
for $(\xi,s)\in{\bf R}^N\times(S,t]$ with $|x|=\epsilon\sqrt{s}$ 
and 
\begin{equation}
\label{eq:5.46}
\tilde{w}(\xi,S)\ge\frac{1}{2}\tilde{\zeta}(S)U(|\xi|)
\ge C^{-1}S^{-\frac{N}{2}}U(\sqrt{S})^{-2}U(|\xi|)
\end{equation}
for $\xi\in B(0,\epsilon\sqrt{S})$. 
By \eqref{eq:5.39}, \eqref{eq:5.44}, \eqref{eq:5.45} and \eqref{eq:5.46}, 
taking a sufficiently large $C_2$ if necessary, 
we see that 
\begin{equation}
\label{eq:5.47}
\tilde{z}\le 0
\end{equation}
for all $(\xi,s)\in{\bf R}^N\times[S,t]$ with $|\xi|=\epsilon\sqrt{s}$ and 
all $(\xi,S)$ with $|\xi|<\epsilon\sqrt{S}$. 
By \eqref{eq:5.43} and \eqref{eq:5.47} 
we apply the comparison principle to obtain 
$\tilde{z}\le 0$ 
for all $(\xi,s)\in{\bf R}^N\times[S,t]$ with $|\xi|\le\epsilon\sqrt{s}$. 
This implies 
$$
p(x,\xi,s)=v(\xi,s)\le C_2U(|x|)\exp\left(-\frac{|x|^2}{C_2t}\right)\tilde{w}(\xi,s)
\le C_2\tilde{\zeta}(s)U(|x|)U(|\xi|)\exp\left(-\frac{|x|^2}{C_2t}\right)
$$
for all $(\xi,s)\in{\bf R}^N\times[S,t]$ with $|\xi|\le\epsilon\sqrt{s}$. 
Taking $(\xi,s)=(y,t)$, by \eqref{eq:5.4} and \eqref{eq:5.42} we obtain
\begin{equation*}
\begin{split}
p(x,y,t) & \le C\tilde{\zeta}(t)U(|x|)U(|y|)\exp\left(-\frac{|x|^2}{Ct}\right)\\
 & \le Ct^{-\frac{N}{2}}\frac{U(\min\{|x|,\sqrt{t}\})U(\min\{|y|,\sqrt{t}\})}{U(\sqrt{t})^2}\exp\left(-\frac{|x-y|^2}{Ct}\right).
\end{split}
\end{equation*}
This means that \eqref{eq:1.6} holds in case~(ii). 
Thus Theorem~\ref{Theorem:1.3} follows.
$\Box$
\section{Positive potentials}\label{section:6}
\noindent
{\bf Proof of Proposition \ref{Proposition:1.1}}. 
We shall use  the classical idea that  
a polynomial decay of a heat kernel is equivalent to a Sobolev inequality. 
We use this to the kernel $G(x,y,t)$ of the operator $L v = - \frac{1}{U^2}\mbox{div}\,(U(x)^2\nabla v)$ 
used  in the proof of Theorem \ref{Theorem:1.2}.  
The $L^p$-spaces in consideration here are $L^p({\bf R}^N, U(x)^2 dx)$ 
and since by assumption $U(x) \sim |x|^\alpha$
a  Sobolev inequality in this setting is a Caffarelli-Kohn-Nirenberg type inequality. 
This strategy was already  used in \cite{Bar} to obtain similar bounds  
for the heat kernel of $-\Delta + \frac{w}{|x|^2}$ for a positive real $w$. 

Let $L$ and $G(x,y,t)$ be as in the proof of Theorem \ref{Theorem:1.2}. Let $\lambda \in {\bf R}$ and  $\phi \in C^\infty({\bf R}^N)$ and bounded with $|\nabla \phi | \le 1$. Let $L_{\lambda, \phi} := e^{-\lambda \phi} L e^{\lambda \phi}$ and $k_{\lambda, \phi}(x,y,t)$ its associated heat kernel. The bilinear form associated to the operator $L_{\lambda, \phi}$ is given by
\begin{eqnarray*}
{\mathcal E}_{\lambda, \phi}(u,v) &=& \int_{{\bf R}^N} (L_{\lambda, \phi} u) v   U(x)^2dx\\
&= & \int_{{\bf R}^N} \nabla(e^{\lambda \phi} u) \nabla (e^{-\lambda \phi} v) U(x)^2 dx\\
&=& \int_{{\bf R}^N} \left[ \nabla u{\cdot}\nabla v + \lambda u \nabla \phi{\cdot}\nabla v 
- \lambda v \nabla \phi{\cdot}\nabla u - \lambda^2 u v \right] U(x)^2 dx.
\end{eqnarray*}
In particular, the quadratic form satisfies  
$$
{\mathcal E}_{\lambda, \phi}(u,u) = \int_{{\bf R}^N} \left[  |\nabla u |^2  - \lambda^2 |u|^2 \right] U(x)^2 dx.
$$
Recall the weighted Sobolev inequality due to  Caffarelli-Korn-Nirenberg \cite{CKN}
\begin{equation}\label{eq:6.1}
\| \nabla u \|^2_{L^2({\bf R}^N, |x|^{2\alpha} dx)} \ge C \| u \|^2_{L^{p_0}({\bf R}^N, |x|^{2\alpha} dx)},
\end{equation}
where $p_0 := \frac{2(N- 2 \alpha)}{N-2-2\alpha}$.\footnote{Here one needs of course $N > 2  + 2 \alpha$. In the case
$N \le 2  + 2 \alpha$, we use a Gagliardo-Nirenberg type inequality instead of \eqref{eq:6.1}. See \cite{CKN}.}
This together with the fact that $U(x) \sim |x|^\alpha$ implies  that
\begin{equation*}
{\mathcal E}_{\lambda, \phi}(u,u) + \lambda^2 \int_{{\bf R}^N}  |u|^2 U(x)^2 dx \ge C \| u \|^2_{L^{p_0}({\bf R}^N, |x|^{2\alpha} dx)}.
\end{equation*}
It is a classical fact that  the semigroup $e^{-t L_{\lambda, \phi}} $ is bounded from $L^2({\bf R}^N, U^2 dx)$ into $L^{p_0}({\bf R}^N, U^2 dx)$ with norm bounded by $C t^{-1/2} e^{\lambda^2t}$. The same strategy as in the proof of a Gaussian upper for the heat kernel of uniformly elliptic operator (see, e.g., \cite{Dav} or \cite{Ouh}) allows to iterate this estimates and see that the semigroup $e^{-t L_{\lambda, \phi}} $ is bounded from $L^2({\bf R}^N, U^2 dx)$ into $L^\infty({\bf R}^N, U^2 dx)$ with norm bounded by $C t^{-N/4 - \alpha/2} e^{\lambda^2t}$. 
Thus,
\begin{equation*}
\int_{{\bf R}^N} | k_{\lambda, \phi}(x,y,t) |^2 U(y)^2 dy \le C t^{-\frac{N}{2} - \alpha} e^{2\lambda^2t}.
\end{equation*}
Set $R_{\lambda, \phi}(x,y,t) := e^{-\lambda \phi(x)}p(x,y,t) e^{\lambda \phi(y)}$. The latter estimate immediately gives
\begin{equation}\label{eq:6.2}
\int_{{\bf R}^N} | R_{\lambda, \phi}(x,y,t) |^2 dy \le C t^{-\frac{N}{2}} \left(\frac{|x|}{\sqrt{t}} \right)^{2\alpha} e^{2\lambda^2t}.
\end{equation}
On the other hand, since $V$ is non-negative we have the domination property
\begin{equation*}
p(x,y,t) \le (4 \pi t)^{-\frac{N}{2}} e^{-\frac{|x-y|^2}{4t}}.
\end{equation*}
This  and the fact that $| \nabla \phi | \le 1$ imply
\begin{equation}\label{eq:6.3}
\int_{{\bf R}^N} | R_{\lambda, \phi}(x,y,t) |^2 dy \le C t^{-\frac{N}{2}} e^{2\lambda^2t}.
\end{equation}
Combining \eqref{eq:6.3} and \eqref{eq:6.2} yields
$$
\int_{{\bf R}^N} | R_{\lambda, \phi}(x,y,t) |^2 dy \le C t^{-\frac{N}{2}} \left(\min(1,\frac{|x|}{\sqrt{t}}) \right)^{2\alpha} e^{2\lambda^2t}.
$$
By the semigroup property and the assumption $U(x) \sim |x|^\alpha$ we have
\begin{eqnarray*}
R_{\lambda,\phi}(x,y,t) &=& \int_{{\bf R}^N}  R_{\lambda, \phi}(x,z,t/2) R_{\lambda, \phi}(z,y,t/2) dz\\
&\le& \left( \int_{{\bf R}^N} | R_{\lambda, \phi}(x,z,t/2) |^2 dz \right)^{1/2} \left( \int_{{\bf R}^N} | R_{\lambda, \phi}(z,y,t/2) |^2 dz \right)^{1/2} \\
&\le&  C t^{-\frac{N}{2}} \left(\min(1,\frac{|x|}{\sqrt{t}}) \right)^{\alpha} \left(\min(1,\frac{|y|}{\sqrt{t}}) \right)^{\alpha} e^{\lambda^2 t}\\
&\le & C t^{-\frac{N}{2}} \frac{U(\min\{|x|,\sqrt{t}\})U(\min\{|y|,\sqrt{t}\})}{U(\sqrt{t})^2} e^{\lambda^2 t}.
\end{eqnarray*}
Hence
\begin{equation*}
p(x,y,t) \le C t^{-\frac{N}{2}}  \frac{U(\min\{|x|,\sqrt{t}\})U(\min\{|y|,\sqrt{t}\})}{U(\sqrt{t})^2} e^{\lambda^2 t}
e^{\lambda(\phi(x) -\phi(y))}.
\end{equation*}
We change $\lambda$ into $-\lambda$ and then optimize as usual over $\lambda$ and $\phi$ to obtain the upper estimate in Theorem \ref{Proposition:1.1}.  
$\Box$\vspace{7pt}

\noindent
{\bf Acknowledgments}. 
The first-named author  was supported partially 
by the Grant-in-Aid for Scientific Research (A)(No.~15H02058) 
from Japan Society for the Promotion of Science. 
The second author was supported partially by JSPS KAKENHI (No.~15K04965 and 15H03631) 
and MEXT KAKENHI (A) (No.~24244012).
The third-named author was  partially supported by the ANR project
HAB, ANR-12-BS01-0013-02.

\bigskip
\noindent Addresses:

\smallskip
\noindent K.I.:  Mathematical Institute, Tohoku University,
Aoba, Sendai 980-8578, Japan.\\
\noindent 
E-mail: {\tt ishige@math.tohoku.ac.jp}\\

\smallskip
\noindent 
Y. K.: Department of Mathematical Sciences, Osaka Prefecture University, 
Sakai 599-8531, Japan. \\
\noindent 
E-mail: {\tt kabeya@ms.osakafu-u.ac.jp}\\

\smallskip
\noindent 
E.M.  O.: Institut de Math\'ematiques, Universit\'e de Bordeaux,  
351, Cours de la Lib\'eration 33405 Talence, France. \\
\noindent 
E-mail: {\tt elmaati.ouhabaz@math.u-bordeaux.fr}\\

\end{document}